\documentclass[12pt, a4paper,reqno]{amsart}
\usepackage[top=25mm, bottom=25mm, left=17mm, right=17mm]{geometry}

\newtheorem{theorem}{Theorem}[section]
\newtheorem{lemma}[theorem]{Lemma}
\newtheorem{proposition}[theorem]{Proposition}
\newtheorem{corollary}[theorem]{Corollary}

\numberwithin{equation}{section}
\numberwithin{table}{section}
\numberwithin{figure}{section}
\newcommand{\abs}[1]{\lvert#1\rvert}
\newcommand{\ve}{\varepsilon}
\newcommand{\ee}{\epsilon}

\newcommand{\wt}[1]{\widetilde{#1}}
\def\ontop#1#2{\setbox0\hbox{#2}\copy0\llap{\raise\ht0\hbox{#1}}}

\newcommand{\vertiii}[1]{{\left\vert\kern-0.25ex\left\vert\kern-0.25ex\left\vert #1 
	\right\vert\kern-0.25ex\right\vert\kern-0.25ex\right\vert}}

\begin{document}
\title{On the well-posedness of the Galerkin semidiscretization of the periodic 
initial-value problem of the Serre equations.}
\author{D.C. Antonopoulos}
\address{Department of Mathematics, University of Athens, 15784 Zographou, Greece, and Institute of 
Applied and Computational Mathematics, FORTH, 70013 Heraklion, Greece}
\email{antonod@math.uoa.gr}
\author{V.A. Dougalis}
\address{Department of Mathematics, University of Athens, 15784 Zographou, Greece, and Institute of 
Applied and Computational Mathematics, FORTH, 70013 Heraklion, Greece}
\email{doug@math.uoa.gr}
\author{D.E. Mitsotakis}
\address{School of Mathematics and Statistics, Victoria University of Wellington, Wellington 6140, 
New Zealand}
\email{dimitrios.mitsotakis@vuw.ac.nz}
\maketitle
\markright{\MakeUppercase{Galerkin semidiscretization of the Serre equations}}
\begin{abstract}
We consider the periodic initial-value problem for the Serre equations of water-wave theory and its 
semidiscrete approximation in the space of smooth periodic polynomial splines. We prove that the 
semidiscrete problem is well posed, locally in time, and satisfies a discrete positivity property 
for the water depth. 
\end{abstract}
\maketitle
\section{Introduction}
In this note we consider the system of Serre equations of water-wave theory, see e.g. \cite{i}, 
\cite{li} and their references, whose initial-value problem is written, for $x\in\mathbb{R}$, 
$t\geq 0$, in unscaled, nondimensional variables, in the form
\begin{equation}
\begin{aligned}
&\zeta_{t} + (\eta u)_{x} = 0,\\
&\eta u_{t} -\tfrac{1}{3}(\eta^{3}u_{tx})_{x} + \eta\zeta_{x} + 
\eta uu_{x} - \tfrac{1}{3}\bigl[\eta^{3}(uu_{xx} - u_{x}^{2})\bigr]_{x} = 0\,,
\end{aligned}
\label{eq11}
\end{equation}
where $\eta = 1+\zeta$, and with given initial data
\begin{equation}
\zeta\big|_{t=0}=\zeta_{0},\quad u\big|_{t=0} = u_{0}.
\label{eq12}
\end{equation}
The solution of this system should satisfy for $t\geq 0$ the positivity condition
\begin{equation}
\exists c_{0}>0, \quad \inf_{x\in\mathbb{R}}\eta(x,t)\geq c_{0}.
\label{eq13}
\end{equation}
Here $x$ and $t$ are proportional to position along an one-dimensional channel and time, 
respectively, $\zeta=\zeta(x,t)$ is the elevation of the free surface at $(x,t)$ above a level of 
rest at height $y=0$, $\eta=1+\zeta$ is the water depth (as, in these variables, the horizontal 
bottom is located at $y=-1$), and $u=u(x,t)$ is the vertically averaged horizontal velocity of the 
water. In \cite{li} Li proved that the Cauchy problem \eqref{eq11}-\eqref{eq12} is well posed, 
locally in time, for $(\zeta,u)\in H^{s}\times H^{s+1}$ for $s>3/2$, where $H^{s}=H^{s}(\mathbb{R})$
is the $L^{2}$-based Sobolev space of functions on $\mathbb{R}$, provided 
$\inf_{x\in\mathbb{R}}\eta_{0}(x)\geq c_{0} >0$; She also proved that the depth positivity condition
\eqref{eq13} is preserved while the solution exists. In \cite{i} Israwi extended Li's result to the 
case of the Serre system with variable bottom topography; for a horizontal bottom the results of 
\cite{li} and \cite{i} coincide. \par
In \cite{adm} we analyzed the Galerkin finite-element semidiscretization of the periodic 
initial-value problem for \eqref{eq11}-\eqref{eq12}. Let $\mathcal{S}_{h}$ be the finite-dimensional
space of smooth, 1-periodic, polynomial splines of order $r\geq 3$ on a uniform mesh in $[0,1]$ with 
meshlength $h$. The semidiscrete approximations $\eta_{h}$, $u_{h}$ of $\eta$, $u$, are 
$\mathcal{S}_{h}$-valued functions of $t\geq 0$ that satisfy for all $\phi$, 
$\chi\in \mathcal{S}_{h}$
\begin{equation}
\begin{aligned}
& (\zeta_{ht},\phi) + \bigl((\eta_{h}u_{h})_{x},\phi\bigr) = 0,\\
& (\eta_{h}u_{ht},\chi) + \tfrac{1}{3}(\eta_{h}^{3}u_{htx},\chi')
+ (\eta_{h}\eta_{hx},\chi) + (\eta_{h}u_{h}u_{hx},\chi) 
+ \tfrac{1}{3}\bigl(\eta_{h}^{3}(u_{h}u_{hxx}-u_{hx}^{2}),\chi'\bigr) = 0,
\end{aligned}
\label{eq14}
\end{equation}
where $\eta_{h}=1+\zeta_{h}$, and $\eta_{h}(0)$, $u_{h}(0)$, are suitable approximations of 
$\eta_{0}$, $u_{0}$ in $\mathcal{S}_{h}$. In \cite{adm} we proved optimal-order $L^{2}$ error 
estimates for the solution of this ode initial-value problem. However, there is a gap in this proof: 
We assumed that this ode initial-value problem has a unique solution $(\eta_{h}(t),u_{h}(t))$
locally in time. In order to ensure this one should prove that $\eta_{h}(t)$ is positive locally in 
time, so that the operator acting on $u_{ht}$ in the second equation of \eqref{eq14} is invertible. 
The purpose of this note is to establish this fact. \par
The proof of local existence of solutions of \eqref{eq11}-\eqref{eq12} and of the positivity 
\eqref{eq13} of the water depth is carried out in \cite{li} and \cite{i} by an analysis of a 
linearized version of the Serre equations and the application of the contraction mapping principle. 
The existence of a solution, as a fixed point of the mapping, is established in \cite{li} using Kato's
general theory, cf. e.g. \cite{k}.
Similarly, a Picard iteration scheme is constructed in \cite{i},
suitable estimates of the iterates are proved, and the convergence of the scheme is established by 
appeal to the `classical method' outlined in \cite{ag} in the context of the proof of local existence of classical
solutions of the initial-value problem for quasilinear symmetric hyperbolic systems.\par
In the paper at hand the proof of local existence and uniqueness for the initial-value problem for 
\eqref{eq14} and of the positivity of $\eta_{h}$, follows, in general lines, the plan of the proof in
\cite{i}. (We also found very helpful as a guide the exposition in Ch. 2 of \cite{m} of the argument
for the local existence of solutions of the initial-value problem for quasilinear hyperbolic systems 
in \cite{k} and \cite{l}.) We construct a sequence $(\eta_{h},u_{h})$ in 
$\mathcal{S}_{h}\times\mathcal{S}_{h}$ satisfying a Picard iteration scheme for the semidiscrete 
problem \eqref{eq14} with the same initial conditions, and prove by induction that, for a small 
enough temporal interval, $(\eta_{h}^{n},u_{h}^{n})$ is bounded in the `high-norm' 
$W^{1,\infty}\times W^{2,\infty}$, while the sequence $(\eta_{ht}^{n},u_{ht}^{n})$ is bounded in 
$L^{\infty}\times W^{1,\infty}$, and $\eta_{h}^{n}$ is bounded below by $c_{0}/8$. The heart of the 
inductive step is an energy proof of the fact that $(\theta^{n},\xi^{n})$ satisfies a high-order
accurate error estimate in $L^{2}\times H^{1}$; here $\theta^{n}=Q_{h}\eta^{n} - \eta_{h}^{n}$, 
$\xi^{n}=Q_{h}u^{n}-u_{h}^{n}$, where $Q_{h}$ is the Thom{\'e}e-Wendroff, \cite{tw}, 
quasiinterpolant in $\mathcal{S}_{h}$, and $(\eta^{n},u^{n})$ is the solution of the analogous 
iteration step set up for pde periodic initial-value problem. (Analogous estimates are satisfied by 
$(\theta_{t}^{n},\xi_{t}^{n})$.) It is finally shown that the sequence
$(\eta_{h}^{n+1}-\eta_{h}^{n},u_{h}^{n+1}-u_{h}^{n})$ is geometrically decreasing in 
$L^{2}\times H^{1}$ and, consequently, that $(\eta_{h}^{n},u_{h}^{n})$ converges as $n\to \infty$ to
the unique solution $(\eta_{h},u_{h})$ of the semidiscrete problem in a small temporal interval.
The solution has the property that $\eta_{h}$ is bounded below by $c_{0}/8$.\par
The contents of the paper are as follows. In section 2 we state the iterative scheme for the 
periodic initial-value problem for the Serre equations and list its properties. In section 3 we 
consider the spline space $\mathcal{S}_{h}$ and present the approximation and stability properties 
of the quasiinterpolant. For a crucial stability estimate that we need, we refer to \cite{dd}. The
basic iterative argument for the local existence of the semidiscrete solution is presented in 
section 4, after same preliminary results that include the proof of consistency of the iterative
scheme. We use the following notation: For integer $k\geq 0$, $C^{k}_{per}=C^{k}_{per}[0,1]$ denotes
the space of 1-periodic functions that are $k$ times continuously differentiable on $[0,1]$, while 
$H_{per}^{k}=H_{per}^{k}(0,1)$ are the $L^{2}$-based Sobolev spaces of 1-periodic functions on 
$[0,1]$; their usual norms are denoted by $\|\cdot\|_{k}$. We let $X_{per}^{k}$ be the space 
$H_{per}^{k}\times H_{per}^{k+1}$ with norm given by 
$\|(v,w)\|_{X_{per}^{k}}=(\|v\|_{k}^{2} + \|w\|_{k+1}^{2})^{1/2}$ for $(v,w)\in X_{per}^{k}$. In 
addition, for $T>0$, $X_{per,T}^{k} = C(0,T;X_{per}^{k})$. The inner product on $L^{2}=L^{2}(0,1)$
is denoted by $(\cdot,\cdot)$ and the associated norm by $\|\cdot\|$. The norms of the spaces 
$W^{k,\infty}$, $L^{\infty}$ will be denoted by $\|\cdot\|_{k,\infty}$, $\|\cdot\|_{\infty}$ 
respectively. Moreover, for real $s>3/2$ we will consider $H^{s}=H^{s}(\mathbb{R})$, the Sobolev
space defined as the subspace of $L^{2}(\mathbb{R})$ consisting of (classes of) functions for which 
$\int_{-\infty}^{\infty}(1+\xi^{2})^{s}\abs{\hat{f}(\xi)}^{2}d\xi<\infty$, where $\hat{f}$ is the 
Fourier transform of $f$; its norm will be denoted again by $\|\cdot\|_{s}$. Finally, 
$H^{s}\times H^{s+1}$ will be denoted by $X^{s}$ with norm 
$\|(v,w)\|_{X^{s}}=(\|v\|_{s}^{2} + \|w\|_{s+1}^{2})^{1/2}$; for $T>0$ we put 
$X_{T}^{s}=C(0,T;X^{s})$.
\section{The initial-value problem} Following the notation of \cite{adm} we write the Cauchy problem
for the Serre equations in nondimensional, unscaled form, for $x\in\mathbb{R}$, $t\geq 0$, as
\begin{align}
& \zeta_{t} + (\eta u)_{x} = 0,
\label{eq21} \\
& \eta u_{t} -\frac{1}{3}(\eta^{3}u_{tx})_{x} + \eta\zeta_{x} + \eta uu_{x}
-\frac{1}{3}\bigl[\eta^{3}(uu_{xx} - u_{x}^{2})\bigr]_{x} = 0,
\label{eq22}
\end{align}
where $\eta = 1 + \zeta$, with given initial values
\begin{equation}
\zeta(x,0)=\zeta_{0}(x),\quad u(x,0)=u_{0}(x),\quad x\in\mathbb{R},
\label{eq23}
\end{equation}
under the hypothesis that the positivity condition
\begin{equation}
\exists\, c_{0}>0, \quad \inf_{x\in\mathbb{R}}\eta \geq c_{0}
\label{eq24}
\end{equation}
holds. As mentioned in the Introduction Li, \cite{li}, proved that the Cauchy problem
\eqref{eq21}-\eqref{eq23} is well posed locally in time for $(\zeta,u)\in X^{s}$, if $s>3/2$, when 
$\inf_{x\in\mathbb{R}}\eta_{0}(x)>0$, and that the property $\inf_{x\in\mathbb{R}}\eta(x,t)>0$ is 
conserved while the solution exists. Subsequently, Israwi, \cite{i}, showed local existence and 
uniqueness of the solution of the Cauchy problem for the Serre system in the case of variable bottom
topography. Since the results of Li and Israwi coincide when the bottom is horizontal, 
\cite[Remark 3]{i} and taking into account \cite[Proposition 1, Appendix B, Theorem 1]{i}, and
\cite[Theorem 5.3]{li}, we may state their results in extended form as in the following theorem.
\begin{theorem} Let $s>3/2$ and $M$ and $c_{0}$ be positive constants such that
\[
\|(\zeta_{0},u_{0})\|X^{s}\leq M/2, \quad \inf_{x\in\mathbb{R}}\eta_{0}(x)\geq c_{0}.
\]
Also let $\zeta^{0}=\zeta_{0}$ $u^{0}=u_{0}$, and for $n=0,1,2,\dots$, consider, for
$x\in\mathbb{R}$, $t\geq 0$, the iterative scheme 
\begin{align}
& \zeta_{t}^{n+1} + \eta^{n}u_{x}^{n+1} + u^{n}\eta_{x}^{n+1} = 0,
\label{eq25}\\
& \eta^{n}u_{t}^{n+1}-\tfrac{1}{3}\bigl((\eta^{n})^{3}u_{tx}^{n+1}\bigr)_{x}
+ \eta^{n}\zeta_{x}^{n+1} + \eta^{n}u^{n}u_{x}^{n+1} 
- \tfrac{1}{3}\bigl[(\eta^{n})^{3}(u^{n}u_{xx}^{n+1}
- u_{x}^{n}u_{x}^{n+1})\bigr]_{x} = 0,
\label{eq26}
\end{align}
with $\eta^{n} = 1 + \zeta^{n}$ and initial values
\begin{equation}
\zeta^{n+1}(x,0)=\zeta_{0}(x), \quad u^{n+1}(x,0) = u_{0}(x).
\label{eq27}
\end{equation}
Then: \\
(i)\,\,
There exists $T_{C}^{0}>0$ such that for $n=0,1,2,\dots$ the Cauchy problem 
\eqref{eq25}-\eqref{eq27} has a unique solution $(\zeta^{n+1},u^{n+1})\in X_{T_{C}^{0}}^{s}$, with 
$(\zeta_{t}^{n+1},u_{t}^{n+1})\in X_{T_{C}^{0}}^{s-1}$, such that
\begin{align*}
& \sup_{0\leq \tau\leq T_{C}^{0}}
\|(\zeta^{n+1}(\tau),u^{n+1}(\tau))\|_{X^{s}}\leq M,\quad
\sup_{0\leq \tau\leq T_{C}^{0}}
\|(\zeta_{t}^{n+1}(\tau),u_{t}^{n+1}(\tau))\|_{X^{s-1}}
\leq M_{1},\\
& \hspace{75pt} \inf_{x\in\mathbb{R}}\eta^{n+1}(x,t) \geq c_{0}/2
\quad \text{for all} \quad t\in [0,T_{C}^{0}],
\end{align*}
where $M_{1}$ is a constant depending on $M$, $c_{0}$.\\
(ii)\,\, There exists $T_{C}>0$, depending on $M$, $c_{0}$, such that the sequence of
solutions 
$(\zeta^{n},u^{n})$ of \eqref{eq25}-\eqref{eq27} converges in $X^{s}$ to the pair
$(\zeta,u)\in X_{T_{C}}^{s}$, with $(\zeta_{t},u_{t})\in X_{T_{C}}^{s-1}$, which is the unique 
solution of the Cauchy problem \eqref{eq21}-\eqref{eq23} in $X_{T}^{s}$, and for which it holds that
\begin{align*}
& \sup_{0\leq \tau\leq T_{C}}
\|(\zeta(\tau),u(\tau))\|_{X^{s}}\leq M,\quad  \sup_{0\leq \tau\leq T_{C}}
\|(\zeta_{t}(\tau),u_{t}(\tau))\|_{X^{s-1}}\leq M_{1}, \\
& \hspace{70pt} \inf_{x\in\mathbb{R}}\eta(x,t) \geq c_{0}/2
\quad \text{for all} \quad t\in [0,T_{C}].
\end{align*}
\end{theorem}
In this paper we consider the following periodic initial-value problem for the Serre equations:
\begin{equation}
\begin{aligned}
& \eta_{t} + (\eta u)_{x} = 0, \\
&\eta u_{t} -\frac{1}{3}\bigl(\eta^{3}u_{xt}\bigr)_{x} + \eta\eta_{x} + \eta uu_{x} 
-\frac{1}{3}\bigl[\eta^{3}(uu_{xx} -u_{x}^{2})\bigr]_{x}= 0,
& \quad x\in [0,1], \quad t\geq 0, \\
& \eta(x,0) =\eta_{0}(x), \quad u(x,0)=u_{0}(x), \quad 0 \leq x\leq 1,  
\end{aligned}
\tag{S}
\label{eqs}
\end{equation}
where $\eta_{0}$, $u_{0}$ are given 1-periodic functions and $\eta_{0}$ satisfies
$\min_{0\leq x\leq 1}\eta_{0}(x)\geq c_{0}>0$ for some constant $c_{0}$. With the aim of proving
local existence and uniqueness of the Galerkin semidiscretization of the problem \eqref{eqs}, we 
assume the validity of the following proposition, which is the periodic analog of the previous
Theorem.
\begin{proposition} Let $s\geq 2$ be an even integer, $(\eta_{0}, u_{0})\in X_{per}^{s}$,
and $M$ and $c_{0}$ positive constants such that
\begin{equation}{\tag{Y$_{0}$}}
\|(\eta_{0},u_{0})\|_{X_{per}^{s}}\leq M/2, \quad \min_{0\leq x\leq 1} \eta_{0}(x) \geq c_{0}.
\label{eqy0}
\end{equation}
Also let $\eta^{0}=\eta_{0}$, $u^{0}=u_{0}$, and for $n=0,1,2,\dots$ consider the iterative scheme
\begin{equation}
\begin{aligned}
& \eta_{t}^{n+1} + \eta^{n}u_{x}^{n+1} + u^{n}\eta_{x}^{n+1} = 0,\\
& \eta^{n}u_{t}^{n+1} - \frac{1}{3}\bigl((\eta^{n})^{3}u_{tx}^{n+1}\bigr)_{x} 
+ \eta^{n}\eta_{x}^{n+1} + \eta^{n}u^{n}u_{x}^{n+1}  
-\frac{1}{3}\bigl[(\eta^{n})^{3}(u^{n}u_{xx}^{n+1} - u_{x}^{n}u_{x}^{n+1})\bigr]_{x} = 0,
\\
& \eta^{n+1}(x,0) =\eta_{0}(x), \quad u^{n+1}(x,0)=u_{0}(x), \, 
\end{aligned}
\tag{S$^{n}$}
\label{eqsn}
\end{equation}
for $x\in [0,1]$, $t\geq 0$. Then: \\
(i)\,\, There exists $T_{per}^{0} > 0$ such that the problem \eqref{eqsn} has a unique solution
$(\eta^{n+1},u^{n+1})$ for $n=0,1,2,\dots$, which belongs to $X_{per,T^{0}_{per}}^{s}$, with
$(\eta_{t}^{n+1},u_{t}^{n+1})\in X_{per,T_{per}^{0}}^{s-1}$, and satisfies
\begin{equation}{\tag{E$^{n}$}}
\begin{aligned}
&\sup_{0\leq\tau\leq T_{per}^{0}}
\|(\eta^{n+1}(\tau),u^{n+1}(\tau))\|_{X_{per}^{s}} \leq M,\\
& \sup_{0\leq \tau\leq T_{per}^{0}}
\|(\eta_{t}^{n+1}(\tau),u_{t}^{n+1}(\tau))\|_{X_{per}^{s-1}}
\leq M_{1},\\
& \min_{0\leq x\leq 1}\eta^{n+1}(x,t) \geq c_{0}/2,
\quad \text{for all} \quad t\in [0,T_{per}^{0}],
\end{aligned}
\label{eqen}
\end{equation}
where $M_{1}$ is a constant that depends on $M$, $c_{0}$.\\
(ii)\,\, There exists $T_{per}>0$, depending on $M$, $c_{0}$, such that the sequence 
$(\eta^{n},u^{n})$ of solutions of \eqref{eqsn} converges in $X_{per}^{s}$ to the pair
$(\eta,u)\in X_{per,T_{per}}^{s}$, with $(\eta_{t},u_{t})\in X_{{per},T_{per}}^{s-1}$ which is the
unique solution of \eqref{eqs} in $X_{per,T_{per}}^{s}$, and which satisfies
\begin{equation}{\tag{E}}
\begin{aligned}
& \sup_{0\leq \tau\leq T_{per}}
\|(\eta(\tau),u(\tau))\|_{X_{per}^{s}}\leq M,\\
& \sup_{0\leq \tau\leq T_{per}}
\|(\eta_{t}(\tau),u_{t}(\tau))\|_{X_{per}^{s-1}}
\leq M_{1}, \\
& \min_{0\leq x\leq 1}\eta(x,t) \geq c_{0}/2
\quad \text{for all} \quad t\in [0,T_{per}].
\end{aligned}
\label{eqe}
\end{equation}
\end{proposition}
\section{Approximation spaces} Let $r$ and $N$ be integers such that $r\geq 3$, $N>4(r-1)$. Let
$h=1/N$ and $x_{i}=ih$, $i=0,1,2,\dots,N$, be a uniform partition of $[0,1]$. We consider the
$N$-dimensional vector space of 1-periodic, smooth, piecewise polynomial splines
\[
\mathcal{S}_{h}=\{v\in C_{per}^{r-2}[0,1] : v_{|(x_{i-1},x_{i})}\in
\mathbb{P}_{r-1}, \,\, i=1,2,\dots,N\},
\]
where $\mathbb{P}_{r-1}$ is the space of polynomials of degree at most $r-1$. It is well known that
$\mathcal{S}_{h}$ has the following approximation properties. Given a sufficiently smooth 1-periodic
function $v$, there exists $\chi\in \mathcal{S}_{h}$ such that
\[
\sum_{j=0}^{q-1}h^{j}\|v-\chi\|_{j} \leq C h^{q} \|v\|_{q}, \quad 1\leq q\leq r,
\]
and
\[
\sum_{j=0}^{q-1}h^{j}\|v-\chi\|_{j,\infty} \leq C h^{q} \|v\|_{q,\infty},
\quad 1\leq q\leq r,
\]
where $C$ is a constant independent of $h$ and $v$. Moreover there exists a constant $C$ 
independent of $h$ such that the inverse properties 
\begin{align*}
\|\chi\|_{\beta} & \leq Ch^{-(\beta - \alpha)} \|\chi\|_{\alpha}, \quad
0\leq \alpha\leq \beta\leq r-1, \\
\|\chi\|_{q,\infty} & \leq Ch^{-(q +1/2)} \|\chi\|, \quad 0\leq q \leq r-1,
\end{align*}
hold for each $\chi\in \mathcal{S}_{h}$. (In the sequel, we will denote by $C$, generically, constants
independent of $h$.) \par
In \cite{tw} Thom\'ee and Wendroff proved that there exists a basis 
$\{\wt{\Phi}_{j}\}_{j=1}^{N}$ of $\mathcal{S}_{h}$ with ${\rm{supp}}(\wt{\Phi}_{j})=O(h)$, such 
that, if $v$ is sufficiently smooth 1-periodic function, the {\em quasiinterpolant}, defined as
$Q_{h}v=\sum_{j=1}^{N}v(x_{j})\wt{\Phi}_{j}$, satisfies 
\begin{equation}
\|Q_{h}v - v\| \leq Ch^{r} \|v^{(r)}\|.
\label{eq31}
\end{equation}
In addition, it was shown in \cite{tw} that the basis $\{\wt{\Phi}_{j}\}_{j=1}^{N}$ may be chosen so
that the following properties hold: \\
(i)\, If $\psi\in \mathcal{S}_{h}$, then
\begin{equation}
\|\psi\| \leq Ch^{-1} \max_{1\leq i\leq N}|(\psi,\wt{\Phi}_{i})|.
\label{eq32}
\end{equation}\indent\,\,
(It follows from \eqref{eq32} that if $\psi\in \mathcal{S}_{h}$, $f\in L^{2}$ are such that 
$(\psi,\wt{\Phi}_{i}) = (f,\wt{\Phi}_{i}) + O(h^{\alpha})$, for 
$1\leq$\\ \indent\,\,\, $i\leq N$, i.e. if $|(\psi - P_{h}f,\wt{\Phi}_{i})| \leq Ch^{\alpha}$, 
$1\leq i\leq N$, where $P_{h}$ is the $L^{2}$-projection operator onto\\ \indent\,\,\, 
$\mathcal{S}_{h}$, then  $\|\psi\|\leq \|\psi - P_{h}f\| + \|P_{h}f\| \leq Ch^{\alpha-1} 
 + \|f\|$.) \\
(ii)\, Let $w$ be a sufficiently smooth 1-periodic function, and $\nu$, $\kappa$ be integers such 
that $0\leq\nu$,\\ \indent\,\,\, $\kappa\leq r-1$. If 
$b_{i} = \bigl( (Q_{h}w)^{(\nu)},\wt{\Phi}_{i}^{(\kappa)}\bigr) 
- (-1)^{\kappa} h w^{(\nu+\kappa)}(x_{i})$,  $1\leq i\leq N$, then
\begin{equation}
\max_{1\leq i\leq N}\abs{b_{i}}
\leq C h^{2r+j-\nu-\kappa}\|w\|_{2r+\nu+\kappa,\infty},
\label{eq33}
\end{equation}\indent \,\,
where $j=1$ if $\nu + \kappa$ is even and $j=2$ if $\nu+\kappa$ is odd. \\
(iii)\, Let $f$, $g$ be sufficiently smooth 1-periodic functions and $\nu$ and $\kappa$ as in (ii) above. 
If \\ \indent\,\,\,
$\beta_{i} = \bigl( f(Q_{h}g)^{(\nu)}, \wt{\Phi}_{i}^{(\kappa)}\bigr) -(-1)^{\kappa}
\bigl(Q_{h}\bigl[(fg^{(\nu)})^{(\kappa)}\bigr],\wt{\Phi}_{i}\bigr)$, $1\leq i\leq N$, then
\begin{equation}
\max_{1\leq i\leq N}|\beta_{i}| \leq C h^{2r+j-\nu-\kappa}
\|f\|_{2r+\kappa,\infty}\|g\|_{2r+\nu+\kappa,\infty},
\label{eq34}
\end{equation}\indent\,\,\,
where $j$ as in (ii). \par
In addition, the following result holds for the 
quasiinterpolant: 
\begin{lemma} Let $r\geq 3$, $v\in H_{per}^{2r}(0,1)\cap W^{2r,\infty}(0,1)$. Then
\begin{align}
& \|Q_{h}v - v\|_{j} \leq Ch^{r-j}\|v\|_{r}, \quad j=0,1,2,
\label{eq35} \\
& \|Q_{h}v - v\|_{j,\infty} \leq Ch^{r-j}\|v\|_{2r,\infty}, \quad j=0,1,2,
\label{eq36}\\
& \|Q_{h}v\|_{j} \leq C\|v\|_{j}, \quad j=1,2, \quad
\|Q_{h}v\|_{j,\infty} \leq C\|v\|_{j,\infty}, \quad j=0,1,2.
\label{eq37}
\end{align}
If in addition $\min_{0\leq x\leq 1}v(x)\geq c_{0}>0$, then there exists $h_{0}$ such that 
\begin{equation}
\min_{0\leq x\leq 1}(Q_{h}v)(x) \geq c_{0}/2, \quad \mbox{for} \quad h
\leq h_{0}.
\label{eq38}
\end{equation}
\end{lemma}
\begin{proof} The estimates \eqref{eq35} follow from the approximation and inverse properties of
$\mathcal{S}_{h}$ and from \eqref{eq31}. The stability estimates \eqref{eq37} follow from 
Proposition 6.1 of \cite{dd}. In order to prove \eqref{eq36} note that from \eqref{eq32}, if $P_{h}$
is the $L^{2}$-projection operator onto $\mathcal{S}_{h}$, it follows for a 1-periodic, sufficiently
smooth function $v$ that
\[
\|P_{h}v - Q_{h}v\|\leq Ch^{-1}\max_{1\leq i\leq N}\abs{(P_{h}v-Q_{h}v,\wt{\Phi}_{i})}
= Ch^{-1}\max_{1\leq i\leq N}\abs{(v-Q_{h}v,\wt{\Phi}_{i})}.
\]
By \eqref{eq34} for $\nu=\kappa=0$, $g=1$, $f=v$, we have therefore 
\[
\|P_{h}v-Q_{h}v\|\leq Ch^{2r}\|v\|_{2r,\infty}.
\]
Hence, for $j=0,1,2$, from the inverse and approximation properties of $\mathcal{S}_{h}$ and the
stability of $P_{h}$ in $W^{j,\infty}$, for $j=0,1,2$ (see \cite[Proposition 5.1]{dd}), it follows
that
\begin{align*}
\|v-Q_{h}v\|_{j,\infty} &\leq \|P_{h}v-Q_{h}v\|_{j,\infty} + \|P_{h}v-v\|_{j,\infty}\\
& \leq Ch^{2r-j-1/2}\|v\|_{2r,\infty} + Ch^{r-j}\|v\|_{r,\infty}\\
&\leq Ch^{r-j}\|v\|_{2r,\infty},
\end{align*}
i.e. that \eqref{eq36} holds. To see \eqref{eq38} note that from \eqref{eq36} it follows that
\[
(Q_{h}v)(x)=(Q_{h}v-v)(x) + v(x)\geq -Ch^{r}\|v\|_{2r,\infty} + c_{0}\geq c_{0}/2,
\]
for $h\leq h_{0}$ with a $h_{0}$ such that $Ch_{0}^{r}\|v\|_{2r,\infty}\leq c_{0}/2$. 
\end{proof} \noindent
{\bf{Remark.}} In the case of piecewise linear, continuous, periodic splines (i.e. $r=2$) it is not
hard to see that regularity of $v$ required in \eqref{eq36} is reduced. Specifically, for 
$v\in C^{2}_{per}$ it holds that
\[
\|v - Q_{h}v\|_{j,\infty}\leq Ch^{2-j}\|v^{(2)}\|_{\infty}, \quad j=0,1.
\]
\section{Galerkin semidiscretization} Let $r\geq 3$ and suppose that for $\eta_{0}$, $u_{0}$
there holds the property \eqref{eqy0} with $s=2r+4$. The Galerkin semidiscretization of \eqref{eqs}
may be defined as follows: For $t\geq 0$ we seek $\eta_{h}$, $u_{h}\in \mathcal{S}_{h}$ such that
\begin{equation}{\tag{S$_{h}$}}
\begin{aligned}
& \eta_{ht} + P_{h}\bigl((\eta_{h}u_{h})_{x}\bigr) = 0,\\
& P_{h}(\eta_{h}u_{ht}) + F_{h}(\eta_{h}^{3}u_{htx}) + P_{h}(\eta_{h}\eta_{hx} 
+ \eta_{h}u_{h}u_{hx}) + F_{h}(\eta_{h}^{3}u_{h}u_{hxx} - \eta_{h}^{3}u_{hx}^{2})=0,\\
& \eta_{h}(0)=Q_{h}\eta_{0},\quad u_{h}(0)=Q_{h}u_{0},
\end{aligned}
\label{eqsh}
\end{equation}
\normalsize
where $P_{h}$ is the $L^{2}$-projection operator onto $\mathcal{S}_{h}$ and
$F_{h} : L^{2}\to \mathcal{S}_{h}$ is defined by
\begin{equation}
(F_{h}(v),\phi) = \frac{1}{3}(v,\phi'), \quad \forall \phi\in \mathcal{S}_{h}.
\label{eq41}
\end{equation}
For the mappings $P_{h}$, $F_{h}$ we will show the following three preliminary Lemmas which will be 
useful in the error estimates in the sequel.
\begin{lemma} Let $w_{h}$, $v_{h}\in \mathcal{S}_{h}$, where in addition $v_{h}\geq c_{0}$ for some
positive constant $c_{0}$. Then
\begin{equation}
\|w_{h}\|_{1} \leq \max\Bigl(\frac{1}{c_{0}},\frac{3}{c_{0}^{3}}\Bigr)
\|P_{h}(v_{h}w_{h}) + F_{h}(v_{h}^{3}w_{hx})\|.
\label{eq42}
\end{equation}
\end{lemma}
\begin{proof} We follow \cite[Lemma 1]{i}. Since
\[
\|w_{h}\|_{1}^{2} \leq \frac{1}{c_{0}}\int_{0}^{1}v_{h}(x)w_{h}^{2}(x)dx
+ \frac{3}{3c_{0}^{3}}\int_{0}^{1}v_{h}^{3}(x)w_{hx}^{2}(x)dx,
\]
there follows that
\begin{equation}
\|w_{h}\|_{1}^{2} \leq \max\Bigl(\frac{1}{c_{0}},\frac{3}{c_{0}^{3}}\Bigr)
(P_{h}(v_{h}w_{h}) + F_{h}(v_{h}^{3}w_{hx}),w_{h}),
\label{eq43}
\end{equation}
from which we get \eqref{eq42}.
\end{proof}
\begin{lemma} Let $v$, $w$, be sufficiently smooth 1-periodic functions and let $\rho=Q_{h}w-w$. Then:\\
$\rm{A}_{1}$. $\|F_{h}(v)\|\leq Ch^{-1}\|v\|$, \\
$\rm{A}_{2}$. $\|P_{h}(v\rho)\| \leq
Ch^{2r}\|v\|_{2r,\infty}\|w\|_{2r,\infty}$,\\
$\rm{A}_{3}$. $\|P_{h}(v\rho_{x})\|\leq
Ch^{2r}\|v\|_{2r,\infty}\|w\|_{2r+1,\infty}$,\\
$\rm{A}_{4}$. $\|F_{h}(v\rho)\|\leq
Ch^{2r}\|v\|_{2r+1,\infty}\|w\|_{2r+1,\infty}$,\\
$\rm{A}_{5}$. $\|F_{h}(v\rho_{x})\|\leq
Ch^{2r-2}\|v\|_{2r+1,\infty}\|w\|_{2r+2,\infty}$,\\
$\rm{A}_{6}$. $\|F_{h}(v\rho_{xx})\|\leq
Ch^{2r-2}\|v\|_{2r+1,\infty}\|w\|_{2r+3,\infty}$
\end{lemma}
\begin{proof} The estimate $\rm{A}_{1}$  follows by the definition \eqref{eq41} of $F_{h}$ and the 
inverse properties of $\mathcal{S}_{h}$.\\
$\rm{A}_{2}$. For $i=1,2,\dots,N$, we have $(v\rho,\wt{\Phi}_{i}) = (vQ_{h}w,\wt{\Phi}_{i}) 
- (vw,\wt{\Phi}_{i}) = (vQ_{h}w,\wt{\Phi}_{i}) - (Q_{h}(vw),\wt{\Phi}_{i})+$\\ \indent\,\,
$(Q_{h}(vw) - vw,\wt{\Phi}_{i})$.
Therefore from \eqref{eq34} we get 
\[
\max_{1\leq i\leq N}\abs{(v\rho,\wt{\Phi}_{i})}
\leq Ch^{2r+1}\|v\|_{2r,\infty}\|w\|_{2r,\infty},
\]\indent\,\,
and \eqref{eq32} yields the desired estimate. \\
$\rm{A}_{3}$. Since $(v\rho_{x},\wt{\Phi}_{i}) = 
(v(Q_{h}w)_{x},\wt{\Phi}_{i}) - (Q_{h}(vw_{x})\wt{\Phi}_{i}) 
+ (Q_{h}(vw_{x}) - vw_{x},\wt{\Phi}_{i})$, for $i=1,2,\dots,N$, \eqref{eq34}\\
\indent\,\,  gives 
$\max_{1\leq i\leq N}\abs{(v\rho_{x},\wt{\Phi}_{i})}
\leq Ch^{2r+1}\|v\|_{2r,\infty}\|w\|_{2r+1,\infty}$; the estimate $\rm{A}_{3}$ 
follows then from \eqref{eq32}. \\ 
$\rm{A}_{4}$. It follows from the fact that 
$F_{h}(v\rho)= - \frac{1}{3}P_{h}(v_{x}\rho + v\rho_{x})$, and from $\rm{A}_{2}$, $\rm{A}_{4}$. \\
$\rm{A}_{5}$. As previously, it suffices to show that 
\[
\max_{1\leq i\leq N}\abs{(v\rho_{x},\wt{\Phi}_{i})} 
\leq Ch^{2r-1}\|v\|_{2r+1,\infty}\|w\|_{2r+2,\infty}.
\] \indent\,\, 
Indeed, since for $i=1,2,\dots,N$,
\[
\quad
(v\rho_{x},\wt{\Phi}_{i}') =
(v(Q_{h}w)_{x},\wt{\Phi}_{i}') + (Q_{h}(vw_{x})_{x},\wt{\Phi}_{i})
-\bigl(Q_{h}(vw_{x})_{x} - (vw_{x})_{x},\wt{\Phi}_{i}\bigr),
\] \indent \,\,
we see that $\rm{A}_{5}$ follows from \eqref{eq34}, \eqref{eq32}, \eqref{eq41}. \\
$\rm{A}_{6}$. For $i=1,2,\dots,N$, we have 
\[ \quad
(v\rho_{xx},\wt{\Phi}_{i}') = (v(Q_{h}w)_{xx},\wt{\Phi}_{i}')
+ (Q_{h}(vw_{xx})_{x},\wt{\Phi}_{i}) -\bigl(Q_{h}(vw_{xx})_{x} - (vw_{xx})_{x},\wt{\Phi}_{i}\bigr).
\]\indent\,\,
Therefore \eqref{eq34}, \eqref{eq32}, \eqref{eq41} yield now $\rm{A}_{6}$.
\end{proof}
A consequence of this Lemma and of \eqref{eq35}-\eqref{eq37} is the following result.
\begin{lemma} Let $u$, $v$, $w$ be 1-periodic, sufficiently smooth functions, and let
$U=Q_{h}u$, $V=Q_{h}v$, $W=Q_{h}w$. Then: \\
$\rm{B}_{1}$. $\|P_{h}(UV - uv)\| \leq
Ch^{2r}\|u\|_{2r,\infty}\|v\|_{2r,\infty}$.\\
$\rm{B}_{2}$. $\|P_{h}(UV_{x}-uv_{x})\| \leq
Ch^{2r-1}\|u\|_{2r,\infty}\|v\|_{2r+1,\infty}$.\\
$\rm{B}_{3}$. $\|P_{h}(UVW_{x}-uvw_{x})\| \leq
Ch^{2r-1}\|u\|_{2r,\infty}\|v\|_{2r,\infty}\|w\|_{2r+1,\infty}$.\\
$\rm{B}_{4}$. $\|F_{h}(U^{3}V_{x} - u^{3}v_{x})\| \leq
Ch^{2r-2}\|u\|_{2r+1,\infty}^{3}\|v\|_{2r+2,\infty}$.\\
$\rm{B}_{5}$. $\abs{\bigl(F_{h}(U^{3}VW_{xx} - u^{3}vw_{xx}),f_{h}\bigr)} \leq C h^{2r-3}
\|u\|_{2r+1,\infty}^{3}\|v\|_{2r+1,\infty}\|w\|_{2r+3,\infty}$.\\
$\rm{B}_{6}$. $\abs{\bigl(F_{h}(U^{3}V_{x}W_{x} - u^{3}v_{x}w_{x}),f_{h}\bigr)}\| \leq Ch^{2r-3}\|u\|_{2r+1,\infty}^{3}\|v\|_{2r+2,\infty}\|w\|_{2r+2,\infty}$.
\end{lemma}
\begin{proof}  Let $\rho^{u} = U-u$, $\rho^{v}=V-v$, and $\rho^{w}=W-w$. Then
\[
UV - uv = (\rho^{u} + u)(\rho^{v} + v) - uv=
\rho^{u}\rho^{v} + \rho^{u}v+ u\rho^{v},
\]
$UV_{x} - uv_{x} = (\rho^{u} + u)(\rho_{x}^{v} + v_{x})-uv_{x}
= \rho^{u}\rho_{x}^{v} + \rho^{u}v_{x} + u\rho_{x}^{v}$,
and, therefore, $\rm{B}_{1}$, $\rm{B}_{2}$ follow from \eqref{eq35}, \eqref{eq36}, and $\rm{A}_{2}$,
$\rm{A}_{3}$. \\
$\rm{B}_{3}$. Here we have 
$UVW_{x} - uvw_{x} = (\rho^{u} + u)(\rho^{v}+v)(\rho_{x}^{w}+w_{x}) - uvw_{x}$, and therefore\\
\indent \,\, 
$UVW_{x} - uvw_{x} = \rho^{u}\rho^{v}W_{x} + \rho^{u}v\rho_{x}^{w}
+ \rho^{u}vw_{x} + u\rho^{v}\rho_{x}^{w}
+u\rho^{v}w_{x} + uv\rho_{x}^{w}$.  Hence 
\begin{align*}
\,\,\,\,\,\, \|P_{h}(UVW_{x} - uvw_{x})\| &\leq
\|\rho^{u}\|_{\infty}\|\rho^{v}\|\|W_{x}\|_{\infty}
+ \|v\|_{\infty}\|\rho^{u}\|_{\infty}\|\rho_{x}^{w}\|
+ \|P_{h}(vw_{x}\rho^{u})\| \\
& \,\,\,\,\,\,
+ \|u\|_{\infty}\|\rho^{v}\|_{\infty}\|\rho_{x}^{w}\|
+ \|P_{h}(uw_{x}\rho^{v})\| + \|P_{h}(uv\rho_{x}^{w})\|,
\end{align*}\indent \,\,
and \eqref{eq35}-\eqref{eq37}, $\rm{A}_{2}$, $\rm{A}_{3}$ yield now the estimate $\rm{B}_{3}$.\\
$\rm{B}_{4}$. It holds that 
\begin{align*}
U^{3}V_{x} - u^{3}v_{x} & =
(\rho^{u} + u)^{3}(\rho_{x}^{v} + v_{x}) - u^{3}v_{x}\\
& = (\rho^{u})^{3}V_{x} + 3(\rho^{u})^{2}uV_{x}
+ 3\rho^{u}u^{2}\rho_{x}^{v}+ 3\rho^{u}u^{2}v_{x} + u^{3}\rho_{x}^{v},
\end{align*}
\indent \,\,
and therefore, taking into account $\rm{A}_{1}$, we get
\begin{align*}
\|F_{h}(U^{3}V_{x} - u^{3}v_{x})\| & \leq
Ch^{-1}(\|(\rho^{u})^{3}\|_{\infty}\|V\|_{1}
+ \|(\rho^{u})^{2}\|_{\infty}\|u\|_{\infty}\|V\|_{1}
+ \|\rho^{u}\|_{\infty}\|u\|_{\infty}^{2}\|\rho^{v}\|_{1})\\
& \,\,\,\,\,\, + 3\|F_{h}(u^{2}v_{x}\rho^{u})\| + \|F_{h}(u^{3}\rho_{x}^{v})\|.
\end{align*}\indent \,\,
This relation and \eqref{eq35}-\eqref{eq37}, and $\rm{A}_{4}$, $\rm{A}_{5}$ give the inequality
$\rm{B}_{4}$. \\
$\rm{B}_{5}$. It holds that 
$U^{3}VW_{xx} - u^{3}vw_{xx} =
(\rho^{u}+u)^{3}(\rho^{v}+v)(\rho_{xx}^{w} + w_{xx})-u^{3}vw_{xx}$, and therefore
\begin{align*}
U^{3}VW_{xx} - u^{3}vw_{xx} & =
(\rho^{u})^{3}VW_{xx} + 3(\rho^{u})^{2}uVW_{xx}
+ 3\rho^{u}u^{2}\rho^{v}W_{xx} \\
& \,\,\,\,\,\, + 3\rho^{u}u^{2}v\rho_{xx}^{w} 
+ 3\rho^{u}u^{2}vw_{xx} + u^{3}\rho^{v}\rho_{xx}^{w} \\
& \,\,\,\,\,\, + u^{3}\rho^{v}w_{xx} + u^{3}v\rho_{xx}^{w}.
\end{align*}\indent \,\,
Hence, in view of $\rm{A}_{1}$, 
\begin{align*}
\qquad \|F_{h}(U^{3}VW_{xx} & - u^{3}vw_{xx})\| \leq Ch^{-1}
(\|\rho^{u}\|_{\infty}^{2}\|\rho^{u}\|\|VW_{xx}\|_{\infty} 
+ \|\rho^{u}\|_{\infty}\|\rho^{u}\|\|u\|_{\infty}\|VW_{xx}\|_{\infty} \\
&+ \|\rho^{u}\|_{\infty}\|u\|_{\infty}^{2}\|\rho^{v}\|\|W_{xx}\|_{\infty}
+ \|\rho^{u}\|_{\infty}\|u\|_{\infty}^{2}\|v\|_{\infty}\|\rho_{xx}^{w}\|
+ \|u\|_{\infty}^{3}\|\rho^{v}\|_{\infty}\|\rho_{xx}^{w}\|)\\
& + 3\|F_{h}(u^{2}vw_{xx}\rho^{u})\| + \|F_{h}(u^{3}w_{xx}\rho^{v})\|
+ \|F_{h}(u^{3}v\rho_{xx}^{w})\|,
\end{align*}\indent\,\,
from which we get $\rm{B}_{5}$ from \eqref{eq35}-\eqref{eq37}, and $\rm{A}_{4}$, $\rm{A}_{6}$.\\
$\rm{B}_{6}$. Here we have
\begin{align*}
\quad U^{3}V_{x}W_{x}-u^{3}v_{x}w_{x} & =
(\rho^{u}+u)^{3}(\rho_{x}^{v} + v_{x})(\rho_{x}^{w}+w_{x})
- u^{3}v_{x}w_{x}\\
& = (\rho^{u})^{3}V_{x}W_{x} + 3(\rho^{u})^{2}uV_{x}W_{x}
+ 3\rho^{u}u^{2}\rho_{x}^{v}W_{x} + 3\rho^{u}u^{2}v_{x}\rho_{x}^{w}\\
& \hspace{10pt}
+ 3\rho^{u}u^{2}v_{x}w_{x} + u^{3}\rho_{x}^{v}\rho_{x}^{w}
+ u^{3}\rho_{x}^{v}w_{x} + u^{3}v_{x}\rho_{x}^{w}.
\end{align*}\indent\,\,
Therefore, as previously, from $\rm{A}_{1}$, \eqref{eq35}-\eqref{eq37}, and $\rm{A}_{4}$, 
$\rm{A}_{5}$, we get $\rm{B}_{6}$.
\end{proof}
In order to prove the existence and uniqueness of solutions of the semidiscrete problem 
\eqref{eqsh}, we define the following iterative scheme, a discrete analog of \eqref{eqsn}:
Let $\eta_{h}^{0}=\eta_{h}(0)$, $u_{h}^{0}=u_{h}(0)$, and for $n=0,1,2,\dots$ consider the system
\begin{equation}{\tag{S$_{h}^{n}$}}
\begin{aligned}
& \eta_{ht}^{n+1}+P_{h}(\eta_{h}^{n}u_{hx}^{n+1})
+ P_{h}(u_{h}^{n}\eta_{hx}^{n+1})=0,\\
& P_{h}(\eta_{h}^{n}u_{ht}^{n+1}) + F_{h}((\eta_{h}^{n})^{3}u_{htx}^{n+1})
+ P_{h}(\eta_{h}^{n}\eta_{hx}^{n+1} + \eta_{h}^{n}u_{h}^{n}u_{hx}^{n+1})
+ F_{h}\bigl((\eta_{h}^{n})^{3}(u_{h}^{n}u_{hxx}^{n+1}-u_{hx}^{n}u_{hx}^{n+1})\bigr)
=0,
\end{aligned}
\label{eqshn}
\end{equation}
with initial values
\[
\eta_{h}^{n+1}(0)=\eta_{h}(0),\quad u_{h}^{n+1}(0)=u_{h}(0).
\]
In the following Lemma we establish a consistency result for the system \eqref{eqshn}, stated in
terms of the quasiinterpolants of the solution $(\eta_{h}^{n},u_{h}^{n})$ of \eqref{eqsn}.
\begin{lemma} Consider $\eta^{n}$, $u^{n}$ as defined by \eqref{eqsn} in Proposition 2.1. Let
$H^{n}=Q_{h}\eta^{n}$, $U^{n}=Q_{h}u^{n}$, and 
$\psi^{n+1}$, $\delta^{n+1} :[0,T_{per}]\to \mathcal{S}_{h}$ be such that 
\begin{equation}
\begin{aligned}
& H_{t}^{n+1} + P_{h}(H^{n}U_{x}^{n+1}) + P_{h}(U^{n}H_{x}^{n+1}) = \psi^{n+1},\\
& P_{h}(H^{n}U_{t}^{n+1})+F_{h}((H^{n})^{3}U_{tx}^{n+1}) + P_{h}(H^{n}H_{x}^{n+1} 
+ H^{n}U^{n}U_{x}^{n+1})\\
& \hspace{63pt}
+ F_{h}\bigl((H^{n})^{3}(U^{n}U_{xx}^{n+1}-U_{x}^{n}U_{x}^{n+1})\bigr) = \delta^{n+1}.
\end{aligned}
\label{eq44}
\end{equation}
Then, there exists a $\wt{h}_{0}>0$, such that 
\begin{equation}
\max_{0\leq t\leq T_{per}}(\|\psi^{n+1}(t)\| + \|\delta^{n+1}(t)\|) \leq C h^{2r-3},
\label{eq45}
\end{equation}
for some constant $C=C(M,M_{1})$, and for $h\leq \wt{h}_{0}$.
\end{lemma}
\begin{proof} Using \eqref{eqy0}, the third inequality of \eqref{eqen} in Proposition 2.1, and
\eqref{eq36}, we see that 
\begin{align*}
& H^{0}=\eta_{0} + Q_{h}\eta_{0} - \eta_{0} \geq c_{0} - Ch^{r}\|\eta_{0}\|_{2r,\infty}, \\
& H^{n} = \eta^{n} + Q_{h}\eta^{n} - \eta^{n} \geq \frac{c_{0}}{2} - Ch^{r}\|\eta^{n}\|_{2r,\infty},
\end{align*}
from which it follows that there exists $\wt{h}_{0}$ such that $H^{0}\geq c_{0}/2$ and
$H^{n}\geq c_{0}/4$ for each $h\leq \wt{h}_{0}$ and for each integer $n\geq 1$.
From the first two equations of \eqref{eqsn} we get
\begin{align*}
& P_{h}\eta_{t}^{n+1} + P_{h}(\eta^{n}u_{x}^{n+1}) + P_{h}(u^{n}\eta_{x}^{n+1}) = 0, \\
& P_{h}(\eta^{n}u_{t}^{n+1}) + F_{h}((\eta^{n})^{3}u_{tx}^{n+1}) 
+ P_{h}(\eta^{n}\eta_{x}^{n+1} + \eta^{n}u^{n}u_{x}^{n+1})
+ F_{h}\bigl((\eta^{n})^{3}(u^{n}u_{xx}^{n+1}-u_{x}^{n}u_{x}^{n+1})\bigr) = 0.
\end{align*}
Subtracting these equations from the analogous ones in \eqref{eq44} yields
\begin{align*}
\psi^{n+1} & = H_{t}^{n+1}-P_{h}\eta_{t}^{n+1}
+ P_{h}(H^{n}U_{x}^{n+1}-\eta^{n}u_{x}^{n+1})
+ P_{h}(U^{n}H_{x}^{n+1} - u^{n}\eta_{x}^{n+1}),\\
\delta^{n+1} & = P_{h}(H^{n}U_{t}^{n+1}-\eta^{n}u_{t}^{n+1})
+ F_{h}((H^{n})^{3}U_{tx}^{n+1}-(\eta^{n})^{3}u_{tx}^{n+1})\\
& \hspace{15pt} + P_{h}(H^{n}H_{x}^{n+1}-\eta^{n}\eta_{x}^{n+1})
+ P_{h}(H^{n}U^{n}U_{x}^{n+1}-\eta^{n}u^{n}u_{x}^{n+1}) \\
& \hspace{15pt} + F_{h}((H^{n})^{3}U^{n}U_{xx}^{n+1}- (\eta^{n})^{3}u^{n}u_{xx}^{n+1})
- F_{h}((H^{n})^{3}U_{x}^{n}U_{x}^{n+1}-(\eta^{n})^{3}u_{x}^{n}u_{x}^{n+1}).
\end{align*}
Therefore, in view of the estimates $\rm{B}_{1}$, $\rm{B}_{2}$ of Lemma 4.3, we see that
\begin{align*}
\|\psi^{n+1}\| \leq Ch^{2r}\|\eta_{t}^{n+1}\|_{2r,\infty}
+ Ch^{2r-1}(& \|\eta^{n}\|_{2r,\infty}\|u^{n+1}\|_{2r+1,\infty}
+ \|u^{n}\|_{2r,\infty}\|\eta^{n+1}\|_{2r+1,\infty}),
\end{align*}
wherefrom, using \eqref{eqen} of Proposition 2.1 with $s=2r+4$, since 
$H_{per}^{s}\subset C_{per}^{s-1}$, we get
\begin{equation}
\|\psi^{n+1}(t)\| \leq C(M,M_{1})h^{2r-1},
\label{eq46}
\end{equation}
for $0\leq t\leq T_{per}$, and for $h\leq \wt{h}_{0}$. In addition, from $\rm{B}_{1}$, $\rm{B}_{4}$,
$\rm{B}_{2}$, $\rm{B}_{3}$, $\rm{B}_{5}$, $\rm{B}_{6}$, of Lemma 4.3 we obtain 
\begin{align*}
\|\delta^{n+1}\| & \leq Ch^{2r}\|\eta^{n}\|_{2r,\infty}
\|u_{t}^{n+1}\|_{2r,\infty}
+ Ch^{2r-2}\|\eta^{n}\|_{2r+1,\infty}^{3}\|u_{t}^{n+1}\|_{2r+2,\infty}\\
& \hspace{15pt} + Ch^{2r-1}\|\eta^{n}\|_{2r,\infty}(\|\eta^{n+1}\|_{2r+1,\infty}
+ \|u^{n}\|_{2r,\infty}\|u^{n+1}\|_{2r+1,\infty})\\
& \hspace{15pt} + Ch^{2r-3}\|\eta^{n}\|^{3}_{2r+1,\infty}\|u^{n}\|_{2r+2,\infty}
\|u^{n+1}\|_{2r+3,\infty},
\end{align*}
and therefore, from \eqref{eqen} it follows that 
\begin{equation}
\|\delta^{n+1}(t)\| \leq C(M,M_{1}) h^{2r-3},
\label{eq47}
\end{equation}
for $0\leq t\leq T_{per}$ and for $h\leq \wt{h}_{0}$. From \eqref{eq47}, \eqref{eq46} we get 
\eqref{eq45}.
\end{proof}
We embark now upon proving the existence-uniqueness of solutions of the iterative scheme 
\eqref{eqshn} and also some several necessary estimates, using induction on $n$. The validity of the
first step of this inductive proof is the content of the following Lemma.
\begin{lemma} Let $\eta^{0}=\eta_{0}$, $u^{0}=u_{0}$, $\eta_{h}^{0}=Q_{h}\eta^{0}$, 
$u_{h}^{0}=Q_{h}u^{0}$, $\eta_{h}^{1}(0)=\eta_{h}(0)$, $u_{h}^{1}(0)=u_{h}(0)$, and $\wt{h}_{0}$ as
in the previous Lemma. Then:\\
(i)\,\, For $h\leq \wt{h}_{0}$ the ode system \eqref{eqshn} for $n=0$ has a unique solution
$(\eta_{h}^{1},u_{h}^{1})$ in  $[0,T_{per}]$. \\
(ii) If $\theta^{1}=Q_{h}\eta^{1}-\eta_{h}^{1}$, $\xi^{1}=Q_{h}u^{1}-u_{h}^{1}$, there exist
constants $\gamma$, $C$, that depend on  $c_{0}$, $M$, $T_{per}$, \\ \indent\,\, such that 
\begin{align}
& \|\theta^{1}\| + \|\xi^{1}\|_{1}\leq \gamma h^{2r-3},
\label{eq48} \\
& \|\theta_{t}^{1}\| + \|\xi_{t}^{1}\|_{1} \leq C h^{2r-4},
\label{eq49}
\end{align}\indent \,\,
for $h\leq \wt{h}_{0}$ and for $t\leq T_{per}$.
\end{lemma}
\begin{proof} (i)\,\, The functions $\eta_{h}^{1}$, $u_{h}^{1}$ satisfy the equations
\begin{equation}
\begin{aligned}
& \eta_{ht}^{1} + P_{h}(H^{0}u_{hx}^{1} + U^{0}\eta_{hx}^{1})=0,\\
& P_{h}(H^{0}u_{ht}^{1}) + F_{h}((H^{0})^{3}u_{htx}^{1})
+ P_{h}(H^{0}\eta_{hx}^{1} + H^{0}U^{0}u_{hx}^{1}) 
+ F_{h}\bigl((H^{0})^{3}(U^{0}u_{hxx}^{1} - U_{x}^{0}u_{hx}^{1})\bigr)=0.
\end{aligned}
\label{eq410}
\end{equation}
Note, as in the beginning of the proof of Lemma 4.4, that 
$H^{0}\geq c_{0}/2$ for $h\leq \wt{h}_{0}$, and therefore the equations of the above linear ode
system are well defined for such values of $h$. Consequently, the initial-value problem for 
\eqref{eq410} has a unique solution $(\eta_{h}^{1},u_{h}^{1})$ for $0\leq t\leq T_{per}$.
Multiplying now the first equation in \eqref{eq410} by $\eta_{h}^{1}$, the second by $u_{h}^{1}$,
and integrating with respect to $x$ on $[0,1]$, gives 
\begin{align*}
& \tfrac{1}{2}\tfrac{d}{dt}\|\eta_{h}^{1}\|^{2}
= - (H^{0}u_{hx}^{1},\eta_{h}^{1}) - (U^{0}\eta_{hx}^{1},\eta_{h}^{1}),\\
& \tfrac{1}{2}\tfrac{d}{dt} \bigl((H^{0}u_{h}^{1},u_{h}^{1})
+ (F_{h}((H^{0})^{3}u_{hx}^{1}),u_{h}^{1})\bigr)
= - (H^{0}\eta_{hx}^{1}+ H^{0}U^{0}u_{hx}^{1},u_{h}^{1})\\
& \hspace{80pt} - \Bigl(F_{h}\bigl((H^{0})^{3}U^{0}u_{hxx}^{1}\bigr),u_{h}^{1}\Bigr)
+ \Bigl(F_{h}\bigl((H^{0})^{3}U_{x}^{0}u_{hx}^{1}\bigr),u_{h}^{1}\Bigr).
\end{align*}
Adding these equations yields 
\begin{align*}
\tfrac{1}{2}\tfrac{d}{dt}\Bigl(\|\eta_{h}^{1}\|^{2} & + (H^{0}u_{h}^{1},u_{h}^{1})  
+ (F_{h}((H^{0})^{3}u_{hx}^{1}),u_{h}^{1})\Bigr) = (H^{0}_{x}u_{h}^{1}, \eta_{h}^{1})
+ \tfrac{1}{2}(U_{x}^{0}\eta_{h}^{1},\eta_{h}^{1})\\
& - (H^{0}U^{0}u_{hx}^{1},u_{h}^{1}) 
- \Bigl(F_{h}\bigl((H^{0})^{3}U^{0}u_{hxx}^{1}\bigr),u_{h}^{1}\Bigr)
+ \Bigl(F_{h}\bigl((H^{0})^{3}U_{x}^{0}u_{hx}^{1}\bigr),u_{h}^{1}\Bigr)\\
& =: a_{1} + a_{2} + a_{3} + a_{4} + a_{5}.
\end{align*}
The terms $a_{i}$ may be suitably estimated by using Lemma 3.1 and the definition of $F_{h}$. We 
first estimate the term $a_{1}+a_{2}+a_{3}$: There holds that
\[
\abs{a_{1}+ a_{2}+ a_{3}} \leq \|H_{x}^{0}\|_{\infty}\|u_{h}^{1}\| \|\eta_{h}^{1}\| +
\|U_{x}^{0}\|_{\infty}\|\eta_{h}^{1}\|^{2} + \|H^{0}\|_{\infty}\|U^{0}\|_{\infty}\|u_{h}^{1}\|^{2},
\]
and therefore, in view of \eqref{eq37}, $\abs{a_{1}+a_{2}+a_{3}} 
\leq C(M) (\|\eta_{h}^{1}\|^{2} + \|u_{h}^{1}\|_{1}^{2})$.\\
In order to estimate $a_{4} + a_{5}$ note that 
\begin{align*}
a_{4} + a_{5} & = - \tfrac{1}{3}((H^{0})^{3}U^{0}u_{hxx}^{1},u_{hx}^{1})
+ \tfrac{1}{3}((H^{0})^{3}U_{x}^{0}u_{hx}^{1},u_{hx}^{1}) \\
& = \tfrac{1}{6}(((H^{0})^{3}U^{0})_{x}u_{hx}^{1},u_{hx}^{1})
+\tfrac{1}{3}((H^{0})^{3}U_{x}^{0}u_{hx}^{1},u_{hx}^{1}),
\end{align*}
and therefore, in view of \eqref{eq37}, $\abs{a_{4}+ a_{5}} \leq C(M) \|u_{h}^{1}\|_{1}^{2}$. These
estimates yield 
\[
\tfrac{1}{2}\tfrac{d}{dt}\Bigl(
\|\eta_{h}^{1}\|^{2}  + (H^{0}u_{h}^{1},u_{h}^{1})
 + (F_{h}((H^{0})^{3}u_{hx}^{1}),u_{h}^{1})\Bigr) \leq C(M)(\|\eta_{h}^{1}\|^{2}
+ \|u_{h}^{1}\|_{1}^{2}),
\]
and, integrating with respect to $t$, we see that 
\begin{align*}
\|\eta_{h}^{1}(t)\|^{2} & + (H^{0}u_{h}^{1},u_{h}^{1}) 
+ \bigl(F_{h}((H^{0})^{3}u_{hx}^{1}),u_{h}^{1}\bigr) \leq \|H^{0}\|^{2} + (H^{0}U^{0},U^{0})\\
& +((F_{h}((H^{0})^{3}U_{x}^{0}),U^{0}) + C(M)\int_{0}^{t}\bigl(\|\eta_{h}^{1}(\tau)\|^{2} 
+ \|u_{h}^{1}(\tau)\|_{1}^{2}\bigr)d\tau.
\end{align*}
Hence, taking into account \eqref{eq42}, we get 
\[
\|\eta_{h}^{1}(t)\|^{2} + \|u_{h}^{1}(t)\|_{1}^{2} \leq C(M,\frac{1}{c_{0}})
\bigl(1 + \int_{0}^{t}\bigl(\|\eta_{h}^{1}(\tau)\|^{2} + 
\|u_{h}^{1}(\tau)\|_{1}^{2}\bigr)d\tau\bigr).
\] 
This inequality and Gronwall's lemma yield the bound 
\[
\|\eta_{h}^{1}(t)\| + \|u_{h}^{1}(t)\|_{1} \leq C(M,\frac{1}{c_{0}},T_{per}),
\quad  t\leq T_{per}.
\]
(ii)\,\, Subtracting the first equation of \eqref{eq410} from the first equation of \eqref{eq44} for
$n=0$, and the second equation of \eqref{eq410} from the second equation of \eqref{eq44} for $n=0$,
we obtain
\begin{equation}
\begin{aligned}
& \theta_{t}^{1} + P_{h}(H^{0}\xi_{x}^{1}) + P_{h}(U^{0}\theta_{x}^{1})
= \psi^{1},\\
& P_{h}(H^{0}\xi_{t}^{1}) + F_{h}((H^{0})^{3}\xi_{tx}^{1})
+ P_{h}(H^{0}\theta_{x}^{1}) + P_{h}(H^{0}U^{0}\xi_{x}^{1})
+ F_{h}((H^{0})^{3}U^{0}\xi_{xx}^{1}) \\
& \hspace{47pt} - F_{h}((H^{0})^{3}U_{x}^{0}\xi_{x}^{1})=\delta^{1}.
\end{aligned}
\label{eq411}
\end{equation}
Multiplying the first equation above by $\theta^{1}$, the second by $\xi^{1}$, and integrating on
$[0,1]$ we get
\begin{align*}
&\tfrac{1}{2}\tfrac{d}{dt}\|\theta^{1}\|^{2} =
-(H^{0}\xi_{x}^{1},\theta^{1}) - (U^{0}\theta_{x}^{1},\theta^{1}) + (\psi^{1},\theta^{1}), \\
& \tfrac{1}{2}\tfrac{d}{dt}\bigl((H^{0}\xi^{1},\xi^{1}) 
+ \bigl(F_{h}((H^{0})^{3}\xi_{x}^{1}),\xi^{1})\bigr) =
-(H^{0}\theta_{x}^{1},\xi^{1}) - (H^{0}U^{0}\xi_{x}^{1},\xi^{1})\\
& \hspace{72pt} - \bigl(F_{h}((H^{0})^{3}U^{0}\xi_{xx}^{1}),\xi^{1}\bigr)
+ \bigl(F_{h}((H^{0})^{3}U_{x}^{0}\xi_{x}^{1}),\xi^{1}\bigr) + (\delta^{1},\xi^{1}).
\end{align*}
Adding these equations yields
\begin{align*}
\tfrac{1}{2}\tfrac{d}{dt}\bigl(\|\theta^{1}\|^{2}
& + (H^{0}\xi^{1},\xi^{1}) + \bigl(F_{h}((H^{0})^{3}\xi_{x}^{1}),\xi^{1}\bigr)
 = (H_{x}^{0}\xi^{1},\theta^{1}) + \tfrac{1}{2}(U_{x}^{0}\theta^{1},\theta^{1})\\
& - \bigl(F_{h}((H^{0})^{3}U^{0}\xi_{xx}^{1}),\xi^{1}\bigr)
+ \bigl(F_{h}((H^{0})^{3}U_{x}^{0}\xi_{x}^{1}),\xi^{1}\bigr) 
+ (\psi^{1},\theta^{1}) + (\delta^{1},\xi^{1}).
\end{align*}
Therefore, using \eqref{eq45}, gives
\[
\tfrac{1}{2}\tfrac{d}{dt}\bigl(\|\theta^{1}\|^{2}
 + (H^{0}\xi^{1},\xi^{1}) + \bigl(F_{h}((H^{0})^{3}\xi_{x}^{1}),\xi^{1})\bigr)
\leq C(M,M_{1})\bigl(\|\theta^{1}\|^{2} + \|\xi^{1}\|_{1}^{2} + h^{4r-6}\bigr),
\]
from which, integrating with respect to $t$ and taking into account \eqref{eq43}, we see that
\[
\|\theta^{1}(t)\|^{2} + \|\xi^{1}(t)\|_{1}^{2} \leq C(M,M_{1},\tfrac{1}{c_{0}})\Bigl(\int_{0}^{t}\bigl(
\|\theta^{1}(\tau)\|^{2} + \|\xi^{1}(\tau)\|_{1}^{2}\bigr)d\tau + h^{4r-6}t\Bigr).
\]
This relation and Gronwall's Lemma implies that there exists a constant $\gamma$ depending on $M$, $M_{1}$,
$c_{0}$, $T_{per}$, such that
\begin{equation}
\max_{0\leq t\leq T_{per}}(\|\theta^{1}(t)\| + \|\xi^{1}(t)\|_{1}) \leq \gamma h^{2r-3},
\label{eq412}
\end{equation}
for $h\leq \wt{h}_{0}$. From \eqref{eq412} we get \eqref{eq48}. Moreover, from the first equation of
\eqref{eq411}, taking into account the inverse properties of $\mathcal{S}_{h}$, \eqref{eq37}, and
\eqref{eq45}, we see that 
\begin{equation}
\max_{0\leq\tau\leq T_{per}}\|\theta_{t}^{1}(\tau)\| \leq C(M,M_{1},\tfrac{1}{c_{0}},T_{per})h^{2r-4},
\label{eq413}
\end{equation}
for $h\leq \wt{h}_{0}$. Multiplication of the second equation of \eqref{eq411} by $\xi_{t}^{1}$ and
integration with respect to $x$ on $[0,1]$ gives
\begin{equation}
\begin{aligned}
\bigl(P_{h}(H^{0}\xi_{t}^{1})& +F_{h}((H^{0})^{3}\xi_{tx}^{1}),
\xi_{t}^{1}\bigr) = (H_{x}^{0}\xi_{t}^{1}+H^{0}\xi_{tx}^{1},\theta^{1})
- (H^{0}U^{0}\xi_{x}^{1},\xi_{t}^{1}) \\
& \hspace{20pt}
+ \bigl(F_{h}((H^{0})^{3}U^{0}\xi_{xx}^{1}),\xi_{t}^{1}\bigr)
+ \bigl(F_{h}((H^{0})^{3}U_{x}^{0}\xi_{x}^{1}),\xi_{t}^{1}\bigr)
+ (\delta^{1},\xi_{t}^{1}).
\end{aligned}
\label{eq414}
\end{equation}
However, 
\begin{align*}
\bigl(F_{h}((H^{0})^{3}U^{0}\xi_{xx}^{1}),\xi_{t}^{1}\bigr)
& = \tfrac{1}{3}\bigl((H^{0})^{3}U^{0}\xi_{xx}^{1},\xi_{tx}^{1}\bigr), \\
\bigl(F_{h}((H^{0})^{3}U_{x}^{0}\xi_{x}^{1}),\xi_{t}^{1}\bigr)
& = \tfrac{1}{3}\bigl((H^{0})^{3}U_{x}^{0}\xi_{x}^{1},\xi_{tx}^{1}\bigr),
\end{align*}
whence from \eqref{eq414} in view of \eqref{eq412}, the inverse properties of $\mathcal{S}_{h}$,
\eqref{eq37}, \eqref{eq45}, and \eqref{eq43} we obtain
\begin{equation}
\max_{0\leq \tau\leq T_{per}}\|\xi_{t}^{1}(\tau)\|_{1} 
\leq C\bigl(M,M_{1},\frac{1}{c_{0}}, T_{per}\bigr)h^{2r-4},
\label{eq415}
\end{equation}
for $h\leq \wt{h}_{0}$. The estimates \eqref{eq413}, \eqref{eq415}, imply \eqref{eq49}.
\end{proof}
In the sequel we will consider the sets 
$\mathcal{A}=\{(v,w)\in \mathcal{S}_{h}\,:\,\|v\|_{1,\infty} + \|w\|_{2,\infty}\leq 3M\}$,
$\mathcal{B}=\{(v,w)\in\mathcal{S}_{h}\,:\, \|v\|_{\infty} + \|w\|_{1,\infty}\leq 3M_{1}\}$, and for
$T>0$, wee will denote $\mathcal{A}_{T}=C(0,T;\mathcal{A})$, $\mathcal{B}_{T}=C(0,T;\mathcal{B})$.
With this notation in place we state the following result, a consequence of the previous Lemma.
\begin{corollary} If $(\eta_{h}^{1},u_{h}^{1})$ is the solution of \eqref{eqshn} for $n=0$, then, 
for $h$ sufficiently small, 
\begin{align}
& (\eta_{h}^{1},u_{h}^{1})\in \mathcal{A}_{T_{per}}, 
\tag{$\Pi_{1}$}
\label{eqp1} \\ 
&(\eta_{ht}^{1},u_{ht}^{1})\in \mathcal{B}_{T_{per}},
\tag{$\Pi_{2}$}  
\label{eqp2} \\
& \min_{0\leq x\leq 1}\eta_{h}^{1}(x,t) \geq c_{0}/8, \quad t\leq T_{per}.
\tag{$\Pi_{3}$}
\label{eqp3}
\end{align}
\end{corollary}
\begin{proof} Taking into account that 
\begin{align*}
\|\eta_{h}^{1}\|_{1,\infty} + \|u_{h}^{1}\|_{2,\infty} & \leq \|\theta^{1}\|_{1,\infty}
+ \|\xi^{1}\|_{2,\infty} + \|Q_{h}\eta^{1}-\eta^{1}\|_{1,\infty}\\
&\,\,\,\,\,\, + \|Q_{h}u^{1}-u^{1}\|_{2,\infty} + \|\eta^{1}\|_{1,\infty} + \|u^{1}\|_{2,\infty},
\end{align*}
and in view of \eqref{eq48}, the approximation properties of $Q_{h}$ in \eqref{eq36}, the inverse
properties of $\mathcal{S}_{h}$, the first inequality of \eqref{eqen}, and the fact that $r\geq 3$,
we see that for sufficiently small $h$, 
\begin{equation}
\max_{0\leq \tau\leq T_{per}}(\|\eta_{h}^{1}(\tau)\|_{1,\infty} +\|u_{h}^{1}(\tau)\|_{2,\infty})
\leq 3M,
\label{eq416}
\end{equation}
that implies \eqref{eqp1}. (For this inequality we need Sobolev's inequality on $[0,1]$.) 
Furthermore, since $\eta_{h}^{1} = H^{1}-\theta^{1}\geq c_{0}/4 - \theta^{1}$, and given that
$\|\theta^{1}(t)\|_{\infty}\leq Ch^{2r-3.5}$ for $t\leq T_{per}$, and $r\geq 3$, we see that
\begin{equation}
\min_{0\leq x\leq 1}\eta_{h}^{1}(x,t) \geq c_{0}/8, 
\label{eq417}
\end{equation}
for sufficiently small $h$ and $t\leq T_{per}$. The inequality \eqref{eq417} is the desired relation
\eqref{eqp3}, from which it also follows that the semidiscrete problem \eqref{eqshn} for $n=1$ has a
unique solution for $h$ sufficiently small. Since
\begin{align*}
\|\eta_{ht}^{1}\|_{\infty} + \|u_{ht}^{1}\|_{1,\infty} & \leq 
\|\theta_{t}^{1}\|_{\infty} + \|\xi_{t}^{1}\|_{1,\infty} 
+ \|Q_{h}\eta_{t}^{1}-\eta_{t}^{1}\|_{\infty}\\
&\,\,\,\,\,\, + \|Q_{h}u_{t}^{1}-u_{t}^{1}\|_{1,\infty} 
+ \|\eta_{t}^{1}\|_{\infty} + \|u_{t}^{1}\|_{1,\infty}, 
\end{align*}
we see from \eqref{eq49}, the approximation properties of $Q_{h}$ in \eqref{eq36}, the inverse
properties of $\mathcal{S}_{h}$, the second inequality of \eqref{eqen}, and the fact that $r\geq 3$,
that
\begin{equation}
\max_{0\leq \tau\leq T_{per}}(\|\eta_{ht}^{1}(\tau)\|_{\infty} +\|u_{ht}^{1}(\tau)\|_{1,\infty})
\leq 3M_{1},
\label{eq418}
\end{equation}
for sufficiently small $h$, i.e. that \eqref{eqp2} holds.
\end{proof}
The following proposition is a basic result of the paper. It establishes, by induction on $n$, the
existence and uniqueness of the solution of \eqref{eqshn}, in an interval $[0,t^{*}]$ for small
enough $t^{*}>0$, and the requisite approximation and boundedness properties of this solution. In 
the sequel we put $\theta^{n} = Q_{h}\eta^{n} - \eta_{h}^{n}$, $\xi^{n} = Q_{h}u^{n}-u_{h}^{n}$.
\begin{proposition} (i)\,\, There exists a $t^{*}$ independent of $h$ such that the problem 
\eqref{eqshn} for $n=1$ has a unique solution $(\eta_{h}^{n+1},u_{h}^{n+1})$, satisfying
$\|\theta^{n+1}(t)\|+\|\xi^{n+1}(t)\|_{1}\leq \gamma h^{2r-3}$ for $t\leq t^{*}$, where $\gamma$ as 
in \eqref{eq48}-\eqref{eq49}, $(\eta_{h}^{n+1},u_{h}^{n+1})\in \mathcal{A}_{t^{*}}$,
$\min_{0\leq x\leq 1}\eta_{h}^{n+1}(x,t)\geq c_{0}/8$, for $0\leq t\leq t^{*}$, and
$(\eta_{ht}^{n+1},u_{ht}^{n+1})\in \mathcal{B}_{t^{*}}$ for sufficiently small $h$.\\
(ii)\,\, For $n\geq 2$: If 
$\|\theta^{n}\| + \|\xi^{n}\|_{1}\leq \gamma h^{2r-3}$, $0\leq t\leq t^{*}$,
$(\eta_{h}^{n},u_{h}^{n})\in \mathcal{A}_{t^{*}}$, 
$\min_{0\leq x\leq 1}\eta_{h}^{n}(x,t)\geq c_{0}/8$, for $0\leq t\leq t^{*}$, and
$(\eta_{ht}^{n},u_{ht}^{n})\in \mathcal{B}_{t^{*}}$, then the problem \eqref{eqshn} has a unique solution 
$(\eta_{h}^{n+1},u_{h}^{n+1})$, such that
\begin{align*}
& \|\theta^{n+1}\| + \|\xi^{n+1}\|_{1} \leq \gamma h^{2r-3}, \quad 0\leq t\leq t^{*},\\
& (\eta_{h}^{n+1},u_{h}^{n+1})\in \mathcal{A}_{t^{*}}, \,\,\,
\min_{0\leq x\leq 1}\eta_{h}^{n}(x,t)\geq c_{0}/8,\,\, 0\leq t\leq t^{*}, \,\,\text{and}\,\,\,
(\eta_{ht}^{n+1},u_{ht}^{n+1})\in \mathcal{B}_{t^{*}}.
\end{align*}
\end{proposition} 
\begin{proof}(i) Let $n=1$. We recall from Lemma 4.5 and Corollary 4.6 that, for sufficiently small $h$,
$\|\theta^{1}\| + \|\xi^{1}\|_{1} \leq \gamma h^{2r-3}$, $0\leq t\leq T_{per}$, and that
$(\eta_{h}^{1},u_{h}^{1})\in \mathcal{A}_{T_{per}}$, 
$(\eta_{ht}^{1},u_{ht}^{1})\in \mathcal{B}_{T_{per}}$, and 
$\min_{0\leq x\leq 1}\eta_{h}^{1}(x,t)\geq c_{0}/8$,\,\, $0\leq t\leq T_{per}$. We also recall that
\eqref{eqp3} of Corollary 4.6 implies that \eqref{eqshn} has a unique solution 
$(\eta_{h}^{n+1},u_{h}^{n+1})$ for $0\leq t\leq T_{per}$. Multiplying the first equation of 
\eqref{eqshn} by $\eta_{h}^{n+1}$, the second by $u_{h}^{n+1}$, using the definition of $F_{h}$ and
integration with respect to $x$ on $[0,1]$, we get 
\begin{equation}
\tfrac{1}{2}\tfrac{d}{dt}\|\eta_{h}^{n+1}\|^{2} = -E_{11} - E_{12},
\label{eq419}
\end{equation}
where $E_{11} = (\eta_{h}^{n}u_{hx}^{n+1},\eta_{h}^{n+1})$, 
$E_{12} = (u_{h}^{n}\eta_{hx}^{n+1},\eta_{h}^{n+1})$, and
\begin{equation}
\tfrac{1}{2}\tfrac{d}{dt}\int_{0}^{1}\bigl[\eta_{h}^{n}(u_{h}^{n+1})^{2} 
+ \tfrac{1}{3}(\eta_{h}^{n})^{3}(u_{hx}^{n+1})^{2}\bigr]dx 
= E_{21} + E_{22} - E_{23} - E_{24} - E_{25} + E_{26},
\label{eq420}
\end{equation}
where $E_{21} = \tfrac{1}{2}(\eta_{ht}^{n}u_{h}^{n+1},u_{h}^{n+1})$,
$E_{22} = \tfrac{1}{2}\bigl(\eta_{ht}^{n}(\eta_{h}^{n})^{2}u_{hx}^{n+1},u_{hx}^{n+1}\bigr)$,
$E_{23} = (\eta_{h}^{n}\eta_{hx}^{n+1},u_{h}^{n+1})$,\\
$E_{24} = (\eta_{h}^{n}u_{h}^{n}u_{hx}^{n+1},u_{h}^{n+1})$,
$E_{25} = \tfrac{1}{3}\bigl((\eta_{h}^{n})^{3}u_{h}^{n}u_{hxx}^{n+1},u_{hx}^{n+1}\bigr)$,
$E_{26} = \tfrac{1}{3}\bigl((\eta_{h}^{n})^{3}u_{hx}^{n}u_{hx}^{n+1},u_{hx}^{n+1}\bigr)$.
Note that $E_{23} = -(\eta_{hx}^{n}\eta_{h}^{n+1},u_{h}^{n+1})-
(\eta_{h}^{n}\eta_{h}^{n+1},u_{hx}^{n+1})$, from which 
$E_{23} =-(\eta_{hx}^{n}\eta_{h}^{n+1},u_{h}^{n+1})- E_{11}$.
In addition
\begin{align*}
E_{25} & = -\tfrac{1}{2}\bigl(\eta_{hx}^{n}(\eta_{h}^{n})^{2}u_{h}^{n}u_{hx}^{n+1},
u_{hx}^{n+1}\bigr) - \tfrac{1}{6}\bigl((\eta_{h}^{n})^{3}u_{hx}^{n}u_{hx}^{n+1},
u_{hx}^{n+1}\bigr) \\
& = -\tfrac{1}{2}\bigl(\eta_{hx}^{n}(\eta_{h}^{n})^{2}u_{h}^{n}u_{hx}^{n+1},
u_{hx}^{n+1}\bigr) - \tfrac{1}{2}E_{26}.
\end{align*}
Adding both sides of \eqref{eq419}, \eqref{eq420} gives
\begin{equation}
\begin{aligned}
\tfrac{1}{2}\tfrac{d}{dt}\Bigl(\|\eta_{h}^{n+1}\|^{2} + \int_{0}^{1}\bigl[
\eta_{h}^{n}(u_{h}^{n+1})^{2} & + \tfrac{1}{3} (\eta_{h}^{n})^{3}(u_{hx}^{n+1})^{2}\bigr]dx\Bigr) 
= E_{21} + E_{22} \\
& + E_{1} - E_{12} - E_{24} + E_{3} + \tfrac{3}{2}E_{26},
\end{aligned}
\label{eq421}
\end{equation}
where $E_{1} = (\eta_{hx}^{n}\eta_{h}^{n+1},u_{h}^{n+1})$, 
$E_{3} = \tfrac{1}{2}\bigl(\eta_{hx}^{n}(\eta_{h}^{n})^{2}u_{h}^{n}u_{hx}^{n+1},u_{hx}^{n+1}\bigr)$.
Since 
\begin{align*}
& \abs{E_{21}} \leq \tfrac{1}{2}\|\eta_{ht}^{n}\|_{\infty}\|u_{h}^{n+1}\|^{2}, \hspace{62pt}
\abs{E_{22}} \leq \tfrac{1}{2}\|\eta_{ht}^{n}\|_{\infty}\|\eta_{h}^{n}\|_{\infty}^{2}
\|u_{hx}^{n+1}\|^{2}, \\
& \abs{E_{1}} \leq \|\eta_{hx}^{n}\|_{\infty}\|\eta_{h}^{n+1}\|\|u_{h}^{n+1}\|, 
\hspace{43pt}
\abs{E_{12}} = \tfrac{1}{2}\abs{(u_{hx}^{n}\eta_{h}^{n+1},\eta_{h}^{n+1})}
\leq \tfrac{1}{2}\|u_{hx}^{n}\|_{\infty}\|\eta_{h}^{n+1}\|^{2},\\
& \abs{E_{24}} \leq \|\eta_{h}^{n}\|_{\infty}\|u_{h}^{n}\|_{\infty}
\|u_{hx}^{n+1}\|\|u_{h}^{n+1}\|, \quad
\abs{E_{3}} \leq \tfrac{1}{2}\|\eta_{hx}^{n}\|_{\infty}\|\eta_{h}^{n}\|_{\infty}^{2}
\|u_{h}^{n}\|_{\infty}\|u_{hx}^{n+1}\|^{2},
\end{align*}
and $\abs{E_{26}} \leq \tfrac{1}{3}\|\eta_{h}^{n}\|_{\infty}^{3}\|u_{hx}^{n}\|_{\infty}
\|u_{hx}^{n+1}\|^{2}$, then, from \eqref{eq421}, Corollary 4.6, and the Cauchy-Schwarz inequality,
it follows that \small
\[
\tfrac{d}{dt}\Bigl(\|\eta_{h}^{n+1}\|^{2} + \int_{0}^{1}\bigl[
\eta_{h}^{n}(u_{h}^{n+1})^{2} + \tfrac{1}{3}
(\eta_{h}^{n})^{3}(u_{hx}^{n+1})^{2}\bigr]dx\Bigr) \leq C(M,M_{1})
(\|\eta_{h}^{n+1}\|^{2} + \|u_{h}^{n+1}\|_{1}^{2}).
\]
\normalsize
Since $\eta_{h}^{n+1}(0)=\eta_{h}^{n}(0)=\eta_{h}(0)$, $u_{h}^{n+1}(0)=u_{h}(0)$,
$\eta_{h}(0)=Q_{h}\eta_{0}$, $u_{h}(0)=Q_{h}u_{0}$, integration with respect to $t$ in the last 
inequality gives, in view of \eqref{eq37},
\begin{align*}
\|\eta_{h}^{n+1}(t)\|^{2} + \int_{0}^{1}\Bigl( \eta_{h}^{n}(u_{h}^{n+1})^{2}
& +\tfrac{1}{3}(\eta_{h}^{n})^{3}(u_{hx}^{n+1})^{2}\Bigr)(x,t)dx
\leq C(M) \\
& + C(M,M_{1})\int_{0}^{t}(\|\eta_{h}^{n+1}(\tau)\|^{2} + \|u_{h}^{n+1}(\tau)\|_{1}^{2})d\tau.
\end{align*}
Therefore, taking into account \eqref{eq42} we see that 
\[
\|\eta_{h}^{n+1}(t)\|^{2} + \|u_{h}^{n+1}(t)\|_{1}^{2} \leq C(M,M_{1},\frac{1}{c_{0}})
\bigl(1 + \int_{0}^{t}\bigl(\|\eta_{h}^{n+1}(\tau)\|^{2} 
+ \|u_{h}^{n+1}(\tau)\|_{1}^{2}\bigr)d\tau\bigr),
\]
from which, using Gronwall's inequality, we get 
\[
\|\eta_{h}^{n+1}(t)\| + \|u_{h}^{n+1}(t)\|_{1} \leq C(M,M_{1},\frac{1}{c_{0}},T_{per}),
\]
for $t\leq T_{per}$. Subtract now the first equation of \eqref{eqshn} from the first equation of
\eqref{eq44} to get
\begin{equation}
\theta_{t}^{n+1} + P_{h}\ve_{11}^{n} + P_{h}\ve_{12}^{n} = \psi^{n+1},
\label{eq422}
\end{equation}
where $\ve_{11}^{n} = H^{n}U_{x}^{n+1} - \eta_{h}^{n}u_{hx}^{n+1}$,
$\ve_{12}^{n} = U^{n}H_{x}^{n+1} - u_{h}^{n}\eta_{hx}^{n+1}$. Subtracting the second equation of
\eqref{eqshn} from the second one of \eqref{eq44} gives
\begin{equation}
P_{h}\ve_{21}^{n} + F_{h}(\ve_{22}^{n}) + P_{h}(\ve_{23}^{n} + \ve_{24}^{n})
+ F_{h}(\ve_{25}^{n}-\ve_{26}^{n}) = \delta^{n+1},
\label{eq423}
\end{equation}
where \small
\begin{equation}
\begin{aligned}
& \ve_{21}^{n} = H^{n}U_{t}^{n+1} - \eta_{h}^{n}u_{ht}^{n+1}, \hspace{63pt}
\ve_{22}^{n} = (H^{n})^{3}U_{tx}^{n+1} - (\eta_{h}^{n})^{3}u_{htx}^{n+1},\\
& \ve_{23}^{n} = H^{n}H_{x}^{n+1} - \eta_{h}^{n}\eta_{hx}^{n+1}, \hspace{63pt}
\ve_{24}^{n} = H^{n}U^{n}U_{x}^{n+1} - \eta_{h}^{n}u_{h}^{n}u_{hx}^{n+1},\\
& \ve_{25}^{n} = (H^{n})^{3}U^{n}U_{xx}^{n+1} - (\eta_{h}^{n})^{3}u_{h}^{n}u_{hxx}^{n+1},\quad
\ve_{26}^{n} = (H^{n})^{3}U_{x}^{n}U_{x}^{n+1} - (\eta_{h}^{n})^{3}u_{hx}^{n}u_{hx}^{n+1}.
\end{aligned}
\label{eq424}
\end{equation}
\normalsize 
Multiply now \eqref{eq422} by $\theta^{n+1}$ and integrate with respect to $x$ on $[0,1]$. This 
gives 
\begin{equation}
\tfrac{1}{2}\tfrac{d}{dt}\|\theta^{n+1}\|^{2} + (\ve_{11}^{n},\theta^{n+1})
+ (\ve_{12}^{n},\theta^{n+1}) = (\psi^{n+1},\theta^{n+1}).
\label{eq425}
\end{equation}
Multiplication of both sides of \eqref{eq423} by $\xi^{n+1}$ and integration with respect to $x$ on
$[0,1]$ yields
\begin{equation}
\begin{aligned}
(\ve_{21}^{n},\xi^{n+1}) + (F_{h}(\ve_{22}),\xi^{n+1}) +
(\ve_{23}^{n}+\ve_{24}^{n},\xi^{n+1})
& + (F_{h}(\ve_{25}^{n}-\ve_{26}^{n}),\xi^{n+1})
 = (\delta^{n+1},\xi^{n+1}).
\end{aligned} 
\label{eq426}
\end{equation}
But $\ve_{21}^{n}=H^{n}U_{t}^{n+1} - \eta_{h}^{n}(U_{t}^{n+1}-\xi_{t}^{n+1})$, and so
$\ve_{21}^{n} = U_{t}^{n+1}\theta^{n} + \eta_{h}^{n}\xi_{t}^{n+1}$, and
\begin{align*}
(\ve_{21}^{n},\xi^{n+1}) & = (U_{t}^{n+1}\theta^{n},\xi^{n+1})
+ (\eta_{h}^{n}\xi_{t}^{n+1},\xi^{n+1})\\
& = (U_{t}^{n+1}\theta^{n},\xi^{n+1}) +
\tfrac{1}{2}\tfrac{d}{dt}(\eta_{h}^{n}\xi^{n+1},\xi^{n+1})
- \tfrac{1}{2}(\eta_{ht}^{n}\xi^{n+1},\xi^{n+1}),
\end{align*}
wherefrom 
\begin{equation}
(\ve_{21}^{n},\xi^{n+1}) = \tfrac{1}{2}\tfrac{d}{dt}(\eta_{h}^{n}\xi^{n+1},\xi^{n+1})
+ (\omega_{1}^{n},\xi^{n+1}),
\label{eq427}
\end{equation}
with
\begin{equation}
\omega_{1}^{n} = U_{t}^{n+1}\theta^{n}-\tfrac{1}{2}\eta_{ht}^{n}\xi^{n+1}.
\label{eq428}
\end{equation}
In addition, $\ve_{22}^{n}=(H^{n})^{3}U_{tx}^{n+1}-(\eta_{h}^{n})^{3} (U_{tx}^{n+1} - 
\xi_{tx}^{n+1})$, i.e. 
$\ve_{22}^{n}=\bigl((H^{n})^{3}-(\eta_{h}^{n})^{3}\bigr)U_{tx}^{n+1} + 
(\eta_{h}^{n})^{3}\xi_{tx}^{n+1}$.
Since, moreover $(H^{n})^{3}-(\eta_{h}^{n})^{3}=\bigl((H^{n})^{2}+H^{n}\eta_{h}^{n}
+(\eta_{h}^{n})^{2}\bigr)\theta^{n}=:a^{n}\theta^{n}$, it follows by the definition of $F_{h}$ in 
\eqref{eq41} that 
$(F_{h}(\ve_{22}^{n}),\xi^{n+1}) =\tfrac{1}{3}(a^{n}U_{tx}^{n+1}\theta^{n},
\xi_{x}^{n+1}) + \tfrac{1}{3}\bigl((\eta_{h}^{n})^{3}\xi_{tx}^{n+1},\xi_{x}^{n+1}\bigr)$,
and therefore
\[
(F_{h}(\ve_{22}^{n}),\xi^{n+1}) =
\tfrac{1}{3}(a^{n}U_{tx}^{n+1}\theta^{n},\xi_{x}^{n+1}) +
\tfrac{1}{6}\tfrac{d}{dt}
\bigl((\eta_{h}^{n})^{3}\xi_{x}^{n+1},\xi_{x}^{n+1}\bigr)
-  \tfrac{1}{2}\bigl(\eta_{ht}^{n}(\eta_{h}^{n})^{2}\xi_{x}^{n+1},\xi_{x}^{n+1}\bigr),
\]
which we write as 
\begin{equation}
(F_{h}(\ve_{22}^{n}),\xi^{n+1}) = \tfrac{1}{6}\tfrac{d}{dt}
\bigl((\eta_{h}^{n})^{3}\xi_{x}^{n+1},\xi_{x}^{n+1}\bigr) + (\omega_{2}^{n},\xi_{x}^{n+1}),
\label{eq429}
\end{equation}
where
\begin{align}
\omega_{2}^{n} & = \tfrac{1}{3}a^{n}
U_{tx}^{n+1}\theta^{n}-\tfrac{1}{2}
\eta_{ht}^{n}(\eta_{h}^{n})^{2}\xi_{x}^{n+1},
\label{eq430}\\
a^{n} & = (H^{n})^{2} + H^{n}\eta_{h}^{n} + (\eta_{h}^{n})^{2}.
\label{eq431}
\end{align}
On the basis of \eqref{eq427} and \eqref{eq429}, \eqref{eq426} may be written as
\begin{align*}
\tfrac{1}{2}\tfrac{d}{dt}\Bigl((\eta_{h}^{n}\xi^{n+1},\xi^{n+1})
& + \tfrac{1}{3}\bigl((\eta_{h}^{n})^{3}\xi_{x}^{n+1},\xi_{x}^{n+1}\bigr)\Bigr)
+ (\omega_{1}^{n},\xi^{n+1}) + (\omega_{2}^{n},\xi_{x}^{n+1})\\
& + (\ve_{23}^{n}+\ve_{24}^{n},\xi^{n+1})
+ (F_{h}(\ve_{25}^{n}-\ve_{26}^{n}),\xi^{n+1}) = (\delta^{n+1},\xi^{n+1}).
\end{align*}
If we add this identity to \eqref{eq425} we finally obtain 
\begin{equation}
\begin{aligned}
\tfrac{1}{2}\tfrac{d}{dt}\Bigl(\|\theta^{n+1}\|^{2} & + 
(\eta_{h}^{n}\xi^{n+1},\xi^{n+1}) +\tfrac{1}{3} 
\bigl((\eta_{h}^{n})^{3}\xi_{x}^{n+1},\xi_{x}^{n+1}\bigr)\Bigr)\\ 
& = -(\ve_{11}^{n}+\ve_{12}^{n},\theta^{n+1}) - (\omega_{1}^{n}+\ve_{23}^{n}+\ve_{24}^{n},\xi^{n+1})
-(\omega_{2}^{n},\xi_{x}^{n+1})\\
&\,\,\,\,\,\,\,  - (F_{h}(\ve_{25}^{n}-\ve_{26}^{n}),\xi^{n+1}) 
+ (\psi^{n+1},\theta^{n+1}) + (\delta^{n+1},\xi^{n+1}).
\end{aligned}
\label{eq432}
\end{equation}
In what follows we proceed to estimate the inner products in the right-hand side of \eqref{eq432},
taking into account the second inequality of \eqref{eq37}, \eqref{eq48}, \eqref{eq49}, and Corollary
4.6. We denote by $\wt{C}$ constants that depend on $c_{0}$, $M$ (and/or $M_{1}$, $\gamma$).\\
$\bullet$\,\, For $(\ve_{11}^{n},\theta^{n+1})$: By the definition of $\ve_{11}^{n}$ after 
\eqref{eq422} we have $\ve_{11}^{n} = H^{n}U_{x}^{n+1} - 
\eta_{h}^{n}(U_{x}^{n+1}-\xi_{x}^{n+1})$, and therefore 
$\ve_{11}^{n} = U_{x}^{n+1}\theta^{n} + \eta_{h}^{n}\xi_{x}^{n+1}$, from which 
$(\ve_{11}^{n},\theta^{n+1}) = (U_{x}^{n+1}\theta^{n},\theta^{n+1})
+ (\eta_{h}^{n}\xi_{x}^{n+1},\theta^{n+1})$. This relation gives
$\abs{(\ve_{11}^{n},\theta^{n+1})} \leq
\|U_{x}^{n+1}\|_{\infty}\|\theta^{n}\|\|\theta^{n+1}\| 
+ \|\eta_{h}^{n}\|_{\infty}\|\|\xi_{x}^{n+1}\|\|\theta^{n+1}\|,
$ i.e. finally
\begin{equation}
\abs{(\ve_{11}^{n},\theta^{n+1})} \leq
\wt{C}(\|\theta^{n+1}\|^{2} + \|\xi^{n+1}\|_{1}^{2} + h^{4r-6}).
\label{eq433}
\end{equation}
$\bullet$\,\, For $(\ve_{12}^{n},\theta^{n+1})$: By the definition of $\ve_{12}^{n}$ in 
\eqref{eq423} we get 
$\ve_{12}^{n} = U^{n}H_{x}^{n+1} - u_{h}^{n}(H_{x}^{n+1}-\theta_{x}^{n+1})$, i.e.
$\ve_{12}^{n} = H_{x}^{n+1}\xi^{n} + u_{h}^{n}\theta_{x}^{n+1}$, and hence
$(\ve_{12}^{n},\theta^{n+1}) = (H_{x}^{n+1}\xi^{n},\theta^{n+1})
+ (u_{h}^{n}\theta_{x}^{n+1},\theta^{n+1})$, i.e. 
\[
(\ve_{12}^{n},\theta^{n+1}) = (H_{x}^{n+1}\xi^{n},\theta^{n+1})
- \tfrac{1}{2} (u_{hx}^{n}\theta^{n+1},\theta^{n+1}).
\] 
This gives 
$\abs{(\ve_{12}^{n},\theta^{n+1})} \leq 
\|H_{x}^{n+1}\|_{\infty}\|\xi^{n}\|\|\theta^{n+1}\| +
\tfrac{1}{2}\|u_{hx}^{n}\|_{\infty}\|\theta^{n+1}\|^{2}$, 
and finally
\begin{equation}
\abs{(\ve_{12}^{n},\theta^{n+1})} \leq \wt{C} (\|\theta^{n+1}\|^{2} + h^{4r-6}).
\label{eq434}
\end{equation}
$\bullet$\,\, For $(\omega_{1}^{n},\xi^{n+1})$: Using the definition of $\omega_{1}^{n}$ in 
\eqref{eq428} we see that
$(\omega_{1}^{n},\xi^{n+1}) = (U_{t}^{n+1}\theta^{n},\xi^{n+1}) -
\tfrac{1}{2}(\eta_{ht}^{n}\xi^{n+1},\xi^{n+1})$, and therefore
\begin{equation}
\abs{(\omega_{1}^{n},\xi^{n+1})} \leq \wt{C}(\|\xi^{n+1}\|^{2} + h^{4r-6}).
\label{eq435}
\end{equation}
$\bullet$\,\, For $(\ve_{23}^{n},\xi^{n+1})$: By the definition of $\ve_{23}^{n}$ in \eqref{eq424}
we have $\ve_{23}^{n}=H^{n}H_{x}^{n+1} - \eta_{h}^{n}(H_{x}^{n+1}-\theta_{x}^{n+1})$, 
i.e. $\ve_{23}^{n}=H_{x}^{n+1}\theta^{n} + \eta_{h}^{n}\theta_{x}^{n+1}$; hence
\begin{align*}
(\ve_{23}^{n},\xi^{n+1}) & = (H_{x}^{n+1}\theta^{n},\xi^{n+1})
+ (\eta_{h}^{n}\theta_{x}^{n+1},\xi^{n+1})
 = (H_{x}^{n+1}\theta^{n},\xi^{n+1}) - (\eta_{hx}^{n}\theta^{n+1},\xi^{n+1})
- (\eta_{h}^{n}\theta^{n+1},\xi_{x}^{n+1}).
\end{align*}
From this we get 
\begin{align*}
\abs{(\ve_{23}^{n},\xi^{n+1})} & \leq
\|H_{x}^{n+1}\|_{\infty}\|\theta^{n}\|\|\xi^{n+1}\|
+ \|\eta_{hx}^{n}\|_{\infty}\|\theta^{n+1}\|\|\xi^{n+1}\|
 + \|\eta_{h}^{n}\|_{\infty}\|\theta^{n+1}\|\|\xi_{x}^{n+1}\|,
\end{align*}
and finally
\begin{equation}
\abs{(\ve_{23}^{n},\xi^{n+1})} \leq \wt{C}
(\|\theta^{n+1}\|^{2} + \|\xi^{n+1}\|_{1}^{2} + h^{4r-6}).
\label{eq436}
\end{equation}
$\bullet$\,\, For $(\ve_{24}^{n},\xi^{n+1})$: From \eqref{eq424} we have
$\ve_{24}^{n}=H^{n}U^{n}U_{x}^{n+1} - \eta_{h}^{n}u_{h}^{n}(U_{x}^{n+1} - \xi_{x}^{n+1})$,
i.e.
$\ve_{24}^{n} =(H^{n}U^{n}-\eta_{h}^{n}u_{h}^{n})U_{x}^{n+1} + \eta_{h}^{n}u_{h}^{n}\xi_{x}^{n+1}$.
But since $H^{n}U^{n}-\eta_{h}^{n}u_{h}^{n}=H^{n}U^{n}- \eta_{h}^{n}(U^{n}-\xi^{n})$,
it follows that 
\[
H^{n}U^{n}-\eta_{h}^{n}u_{h}^{n}=U^{n}\theta^{n}+\eta_{h}^{n}\xi^{n},
\]
and
$\ve_{24}^{n}=U^{n}U_{x}^{n+1}\theta^{n}+\eta_{h}^{n}U_{x}^{n+1}\xi^{n}
+ \eta_{h}^{n}u_{h}^{n}\xi_{x}^{n+1}$.
Therefore 
\[
(\ve_{24}^{n},\xi^{n+1}) =(U^{n}U_{x}^{n+1}\theta^{n},\xi^{n+1})
+ (\eta_{h}^{n}U_{x}^{n+1}\xi^{n},\xi^{n+1})
+ (\eta_{h}^{n}u_{h}^{n}\xi_{x}^{n+1},\xi^{n+1}),
\]
from which
\[
\abs{(\ve_{24}^{n},\xi^{n+1})} \leq
\bigl(\|U^{n}\|_{\infty}\|U_{x}^{n+1}\|_{\infty}\|\theta^{n}\| +
\|\eta_{h}^{n}\|_{\infty}\|U_{x}^{n+1}\|_{\infty}\|\xi^{n}\| 
+ \|\eta_{h}^{n}\|_{\infty}\|u_{h}^{n}\|_{\infty}\|\xi_{x}^{n+1}\|\bigr)\|\xi^{n+1}\|.
\]
Hence
\begin{equation}
\abs{(\ve_{24}^{n},\xi^{n+1})} \leq \wt{C}(\|\xi^{n+1}\|_{1}^{2} + h^{4r-6}).
\label{eq437}
\end{equation}
$\bullet$\,\, For $(\omega_{2}^{n},\xi_{x}^{n+1})$: Using the definition of $\omega_{2}^{n}$ in
\eqref{eq430}, gives
\begin{align*}
\abs{(\omega_{2}^{n},\xi_{x}^{n+1})} & \leq \tfrac{1}{3}\bigl(
\|H^{n}\|_{\infty}^{2} + \|H^{n}\|_{\infty}\|\eta_{h}^{n}\|_{\infty}
+ \|\eta_{h}^{n}\|^{2}_{\infty}\bigr)\|U_{tx}^{n}\|_{\infty}\|\theta^{n}\|
\|\xi_{x}^{n+1}\| \\
& \,\,\,\,\,\,+ \tfrac{1}{3}\|\eta_{ht}^{n}\|_{\infty}\|\eta_{h}^{n}\|^{2}_{\infty}
\|\xi_{x}^{n+1}\|^{2},
\end{align*}
i.e. that
\begin{equation}
\abs{(\omega_{2}^{n},\xi_{x}^{n+1})} \leq \wt{C} (\|\xi^{n+1}\|_{1}^{2} + h^{4r-6}).
\label{438}
\end{equation}
$\bullet$\,\, For $(F_{h}(\ve_{25}^{n}),\xi^{n+1})$: From the definition of $\ve_{25}^{n}$ in 
\eqref{eq424} we see that
\begin{align*}
\ve_{25}^{n} & = (H^{n})^{3}U^{n}U_{xx}^{n+1} -
(\eta_{h}^{n})^{3}u_{h}^{n}(U_{xx}^{n+1} - \xi_{xx}^{n+1})\\
& = \bigl((H^{n})^{3}U^{n} -(\eta_{h}^{n})^{3}u_{h}^{n}\bigr)U_{xx}^{n+1}
+ (\eta_{h}^{n})^{3}u_{h}^{n}\xi_{xx}^{n+1},
\end{align*} 
and since 
\begin{align*}
(H^{n})^{3}U^{n} - (\eta_{h}^{n})^{3}u_{h}^{n} & =
(H^{n})^{3}U^{n} - (\eta_{h}^{n})^{3}(U^{n}-\xi^{n})\\
& = \bigl((H^{n})^{3}-(\eta_{h}^{n})^{3}\bigr)U^{n}
+ (\eta_{h}^{n})^{3}\xi^{n}\\
& = \bigl((H^{n})^{2} + H^{n}\eta_{h}^{n} + (\eta_{h}^{n})^{2}\bigr)
U^{n}\theta^{n} + (\eta_{h}^{n})^{3}\xi^{n},
\end{align*} 
we get $\ve_{25}^{n} = a^{n}U^{n}U_{xx}^{n+1}\theta^{n} + (\eta_{h}^{n})^{3}U_{xx}^{n+1}\xi^{n}
+ (\eta_{h}^{n})^{3}u_{h}^{n}\xi_{xx}^{n+1}$, with $a^{n}$ as in \eqref{eq431}. 
Hence, since $(F_{h}(\ve_{25}^{n}),\xi^{n+1})=\tfrac{1}{3}(\ve_{25}^{n},\xi_{x}^{n+1})$, we obtain 
\begin{align*}
(F_{h}(\ve_{25}^{n}),\xi^{n+1}) & = \tfrac{1}{3}(a^{n}U^{n}U_{xx}^{n+1}\theta^{n},\xi_{x}^{n+1})
+ \tfrac{1}{3}\bigl((\eta_{h}^{n})^{3}U_{xx}^{n+1}\xi^{n},\xi_{x}^{n+1}\bigr)\\
&\,\,\,\,\,\, + \tfrac{1}{3} \bigl((\eta_{h}^{n})^{3}u_{h}^{n}\xi_{xx}^{n+1},\xi_{x}^{n+1}\bigr)\\
& = \tfrac{1}{3}(a^{n}U^{n}U_{xx}^{n+1}\theta^{n},\xi_{x}^{n+1})
+ \tfrac{1}{3}\bigl((\eta_{h}^{n})^{3}U_{xx}^{n+1}\xi^{n},\xi_{x}^{n+1}\bigr)\\
& \,\,\,\,\,\, - \tfrac{1}{2}
\bigl(\eta_{hx}^{n}(\eta_{h}^{n})^{2}u_{h}^{n}\xi_{x}^{n+1},\xi_{x}^{n+1}\bigr)
-\tfrac{1}{6}
\bigl((\eta_{h}^{n})^{3}u_{hx}^{n}\xi_{x}^{n+1},\xi_{x}^{n+1}\bigr).
\end{align*} 
It follows that 
\begin{equation}
\abs{(F_{h}(\ve_{25}^{n}),\xi^{n+1})} \leq \wt{C}(\|\xi^{n+1}\|_{1}^{2}
+ h^{4r-6}).
\label{eq439}
\end{equation}
$\bullet$\,\, For $(F_{h}(\ve_{26}^{n}),\xi^{n+1})$: Using the definition of $\ve_{26}^{n}$ in 
\eqref{eq424} we get
\begin{align*}
\ve_{26}^{n} & = (H^{n})^{3}U_{x}^{n}U_{x}^{n+1} -
(\eta_{h}^{n})^{3}u_{hx}^{n}(U_{x}^{n+1} - \xi_{x}^{n+1})\\
& = \bigl((H^{n})^{3}U_{x}^{n} - (\eta_{h}^{n})^{3}u_{hx}^{n}\bigr)U_{x}^{n+1}
+ (\eta_{h}^{n})^{3}u_{hx}^{n}\xi_{x}^{n+1},
\end{align*} 
and since 
\begin{align*}
(H^{n})^{3}U_{x}^{n} - (\eta_{h}^{n})^{3}u_{hx}^{n} & = (H^{n})^{3}U_{x}^{n} 
- (\eta_{h}^{n})^{3}(U_{x}^{n}-\xi_{x}^{n})\\
& =\bigl((H^{n})^{3}-(\eta_{h}^{n})^{3}\bigr)U_{x}^{n} + (\eta_{h}^{n})^{3}\xi_{x}^{n},
\end{align*}
from which $(H^{n})^{3}U_{x}^{n} - (\eta_{h}^{n})^{3}u_{hx}^{n} = a^{n}U_{x}^{n}\theta^{n} 
+ (\eta_{h}^{n})^{3}\xi_{x}^{n}$,
we get 
$\ve_{26}^{n} = 
a^{n}U_{x}^{n}U_{x}^{n+1}\theta^{n} + (\eta_{h}^{n})^{3}U_{x}^{n+1}\xi_{x}^{n}
+ (\eta_{h}^{n})^{3}u_{hx}^{n}\xi_{x}^{n+1}$. Hence
\begin{align*}
(F_{h}(\ve_{26}^{n}),\xi^{n+1}) & = \tfrac{1}{3}(a^{n}U_{x}^{n}U_{x}^{n+1}\theta^{n},\xi_{x}^{n+1})
+ \tfrac{1}{3}\bigl((\eta_{h})^{3}U_{x}^{n+1}\xi_{x}^{n},\xi_{x}^{n+1}\bigr)
+ \tfrac{1}{3} \bigl((\eta_{h}^{n})^{3}u_{hx}^{n}\xi_{x}^{n+1},\xi_{x}^{n+1}\bigr).
\end{align*} 
Thus 
\begin{equation}
\abs{(F_{h}(\ve_{26}^{n}),\xi^{n+1})} \leq \wt{C} (\|\xi^{n+1}\|_{1}^{2} + h^{4r-6}).
\label{eq440}
\end{equation}
From \eqref{eq432}, taking into account \eqref{eq433}-\eqref{eq440}, and \eqref{eq45}, we get
\[
\tfrac{d}{dt}\Bigl(\|\theta^{n+1}\|^{2} + (\eta_{h}^{n}\xi^{n+1},\xi^{n+1}) +
\tfrac{1}{3}\bigl((\eta_{h}^{n})^{3}\xi_{x}^{n+1},\xi_{x}^{n+1}\bigr)\Bigr)
\leq \wt{C}\bigl(\|\theta^{n+1}\|^{2} + \|\xi^{n+1}\|_{1}^{2} + h^{4r-6}\bigr).
\]
Integration with respect to $t$ in the above gives, in view of \eqref{eq42}, that
\[
\|\theta^{n+1}(t)\|^{2} + \|\xi^{n+1}(t)\|_{1}^{2} 
\leq \wt{C}\Bigl(\int_{0}^{t}\bigl(\|\theta^{n+1}(\tau)\|^{2}
+ \|\xi^{n+1}(\tau)\|_{1}^{2}\bigr)d\tau + h^{4r-6}t\Bigr).
\]
From this relation and Gronwall's Lemma it follows that 
\[
\|\theta^{n+1}(t)\| + \|\xi^{n+1}(t)\|_{1} \leq \sqrt{2}h^{2r-3}(e^{\wt{C}t} - 1)^{1/2},
\]
for $t\leq T_{per}$. Therefore, if we choose $t^{*}$ such that 
$\sqrt{2}(e^{\wt{C}t^{*}}-1)^{1/2}\leq \gamma$, where $\gamma$ as in \eqref{eq48}-\eqref{eq49}, then
\begin{equation}
\|\theta^{n+1}(t)\| + \|\xi^{n+1}(t)\|_{1} \leq \gamma h^{2r-3},
\label{eq441}  
\end{equation}
for $t\leq t^{*}$. \par
In order to complete (i) of Proposition 4.7 we must show that (for $n=1$)
$(\eta_{h}^{n+1},u_{h}^{n+1})\in \mathcal{A}_{t^{*}}$. This follows from \eqref{eq441}
if we argue as in the proof of \eqref{eqp1} in Corollary 4.6. The fact that
$\min_{0\leq x\leq 1}\eta_{h}^{n+1}(x,t)\geq c_{0}/8$, $t\leq t^{*}$, follows from \eqref{eq441} if
we argue as in the proof of \eqref{eqp3} in the same Corollary. It remains to show that
$(\eta_{ht}^{n+1},u_{ht}^{n+1})\in \mathcal{B}_{t^{*}}$. To see this, note that from \eqref{eq422}
the fact that  $\ve_{11}^{n}=U_{x}^{n+1}\theta^{n} + \eta_{h}^{n}\xi_{x}^{n+1}$,
$\ve_{12}^{n}=H_{x}^{n+1}\xi^{n} + u_{h}^{n}\theta_{x}^{n+1}$, \eqref{eq441}, the inverse properties
of $\mathcal{S}_{h}$, \eqref{eq45}, Corollary 4.6, and the second inequality of \eqref{eq37}, it
follows that there exists a constant $C$ such that for $0\leq t\leq t^{*}$
\begin{equation}
\|\theta_{t}^{n+1}\|\leq C h^{2r-4}.
\label{eq442}
\end{equation}
Taking $L^{2}$-inner products of both sides of \eqref{eq423} with $\xi_{t}^{n+1}$ and in view of the
definition of $F_{h}$ in \eqref{eq41}, we see that
\begin{equation}
(\ve_{21}^{n},\xi_{t}^{n+1}) + \tfrac{1}{3}(\ve_{22}^{n},\xi_{tx}^{n+1})
+ (\ve_{23} + \ve_{24},\xi_{t}^{n+1})
+ \tfrac{1}{3}(\ve_{25}^{n}-\ve_{26}^{n},\xi_{tx}^{n+1})
=(\delta^{n+1},\xi_{t}^{n+1}).
\label{eq443}
\end{equation}
Since $\ve_{21}^{n}=U_{t}^{n+1}\theta^{n} + \eta_{h}^{n}\xi_{t}^{n+1}$,
$\ve_{22}^{n} = a^{n}U_{tx}^{n+1}\theta^{n} + (\eta_{h}^{n})^{3}\xi_{tx}^{n+1}$, we may write
\eqref{eq443} as
\begin{equation}
\begin{aligned}
(\eta_{h}^{n}\xi_{t}^{n+1},\xi_{t}^{n+1}) & + 
\tfrac{1}{3}\bigl((\eta_{h}^{n})^{3}\xi_{tx}^{n+1},\xi_{tx}^{n+1}\bigr)
= - (U_{t}^{n+1}\theta^{n},\xi_{t}^{n+1}) - 
\tfrac{1}{3}(a^{n}U_{tx}^{n+1}\theta^{n},\xi_{tx}^{n+1})\\
& - (\ve_{23}^{n} + \ve_{24}^{n},\xi_{t}^{n+1}) 
- \tfrac{1}{3}(\ve_{25}^{n}-\ve_{26}^{n},\xi_{tx}^{n+1}) + (\delta^{n+1},\xi_{t}^{n+1}).
\end{aligned}
\label{eq444}
\end{equation} 
Since $\ve_{23}^{n} = H_{x}^{n+1}\theta^{n} + \eta_{h}^{n}\theta_{x}^{n+1}$, 
$\ve_{24}^{n} = U^{n}U_{x}^{n+1}\theta^{n} + \eta_{h}^{n}U_{x}^{n+1}\xi^{n}
+\eta_{h}^{n}u_{h}^{n}\xi_{x}^{n+1}$, and 
\begin{align*}
\ve_{25}^{n} & = a^{n}U^{n}U_{xx}^{n+1}\theta^{n} + (\eta_{h}^{n})^{3}U_{xx}^{n+1}\xi^{n}
+ (\eta_{h}^{n})^{3}u_{h}^{n}\xi_{xx}^{n+1},\\
\ve_{26}^{n} & = a^{n}U_{x}^{n}U_{x}^{n+1}\theta^{n} + (\eta_{h}^{n})^{3}U_{x}^{n+1}\xi_{x}^{n}
+ (\eta_{h}^{n})^{3}u_{hx}^{n}\xi_{x}^{n+1},
\end{align*}
it follows from \eqref{eq444}, \eqref{eq441}, Corollary 4.6, the inverse properties of
$\mathcal{S}_{h}$, and \eqref{eq37}, that there exists a constant $C$ such that for 
$0\leq t\leq t^{*}$
\begin{equation}
\|\xi_{t}^{n+1}\|_{1} \leq Ch^{2r-4}.
\label{eq445}
\end{equation}
From \eqref{eq442}, \eqref{eq445}, and arguments similar to those as of the proof of Corollary 4.6,
it follows that $(\eta_{ht}^{n+1},u_{ht}^{n+1})\in \mathcal{B}_{t^{*}}$. We conclude that (i) of 
Proposition 4.7 is valid. \\
(ii) The proof of the inductive step for $n\geq 2$ follows  mutatis mutandis that of part (i) if we
replace $T_{per}$ with $t^{*}$.
\end{proof}
On the basis of the results of Lemma 4.5, Corollary 4.6, and Proposition 4.7 we may state the 
following
\begin{corollary} Let $h$ be sufficiently small. Then there exists $t^{*}>0$, 
and a constant $C$, both independent of $h$ and $n$, such that the solutions $(\eta_{h}^{n},u_{h}^{n})$ of \eqref{eqshn} for 
$n=1,2,\dots$ exist uniquely for $0\leq t\leq t^{*}$ and satisfy\\
{\rm{(i)}}\,\,\,\,\, $\|\theta^{n}\| + \|\xi^{n}\|_{1} \leq \gamma h^{2r-3}, 
\quad 0\leq t\leq t^{*}$, \\
{\rm{(ii)}}\,\,\,\, $\|\theta_{t}^{n}\| + \|\xi^{n}_{t}\|_{1} \leq C h^{2r-4}, 
\quad 0\leq t\leq t^{*}$,\\
{\rm{(iii)}}\,\, $(\eta_{h}^{n},u_{h}^{n})\in \mathcal{A}_{t^{*}}$, 
$(\eta_{ht}^{n},u_{ht}^{n})\in \mathcal{B}_{t^{*}}$,\\
{\rm{(iv)}}\,\,\, $\min_{0\leq x\leq 1}\eta_{h}^{n}(x,t)\geq c_{0}/8, \quad 0\leq t\leq t^{*}$.
\end{corollary}
The final step of the proof consists in showing that the sequence $(\eta_{h}^{n},u_{h}^{n})$
converges, as $n\to \infty$, to the solution $(\eta_{h},u_{h})$ of \eqref{eqsh}. We first establish 
that the sequences $(\eta_{h}^{n+1}-\eta_{h}^{n},u_{h}^{n+1}-u_{h}^{n})$ diminish geometrically
in a suitable norm.
\begin{proposition} Let $\delta\eta_{h}^{n}:=\eta_{h}^{n+1}-\eta_{h}^{n}$, 
$\delta u_{h}^{n}:=u_{h}^{n+1} - u_{h}^{n}$. Then there exist $t_{1}^{*}>0$, and  $\alpha\in (0,1)$, both independent of $h$ and $n$, 
such that for all $n\geq 1$
\begin{equation}
\max_{0\leq t\leq t_{1}^{*}}(\|\delta \eta_{h}^{n}(t)\|^{2} + \|\delta u_{h}^{n}(t)\|_{1}^{2})^{1/2} 
\leq \alpha \max_{0\leq t\leq t_{1}^{*}}(\|\delta \eta_{h}^{n-1}(t)\|^{2} + 
\|\delta u_{h}^{n-1}(t)\|_{1}^{2})^{1/2}.
\label{eq446}
\end{equation}
\end{proposition}
\begin{proof} Subtracting the equations of (S$_{h}^{n-1}$) from the analogous equations of \eqref{eqshn} we
have
\begin{equation}
\begin{aligned}
& \delta\eta_{ht}^{n} = -P_{h}\epsilon_{11}^{n} - P_{h}\epsilon_{12}^{n}, \\
& P_{h}\epsilon_{21}^{n}= - F_{h}(\epsilon_{22}^{n}) - P_{h}(\epsilon_{23}^{n}+\epsilon_{24}^{n})
- F_{h}(\epsilon_{25}^{n}-\epsilon_{26}^{n}),
\end{aligned}
\label{eq447}
\end{equation}
where 
\begin{align*}
\ee_{11}^{n} & = \eta_{h}^{n}u_{hx}^{n+1} - \eta_{h}^{n-1}u_{hx}^{n},\hspace{75pt}
\ee_{12}^{n} = u_{h}^{n}\eta_{hx}^{n+1} - u_{h}^{n-1}\eta_{hx}^{n}, \\
\ee_{21}^{n} & = \eta_{h}^{n}u_{ht}^{n+1} - \eta_{h}^{n-1}u_{ht}^{n},\hspace{75pt}
\ee_{22}^{n} = (\eta_{h}^{n})^{3}u_{htx}^{n+1} - (\eta_{h}^{n-1})^{3}u_{htx}^{n},\\
\ee_{23}^{n} & = \eta_{h}^{n}\eta_{hx}^{n+1} - \eta_{h}^{n-1}\eta_{hx}^{n},\hspace{75pt}
\ee_{24}^{n} = \eta_{h}^{n}u_{h}^{n}u_{hx}^{n+1} - \eta_{h}^{n-1}u_{h}^{n-1}u_{hx}^{n},\\
\ee_{25}^{n} & = (\eta_{h}^{n})^{3}u_{h}^{n}u_{hxx}^{n+1}
-(\eta_{h}^{n-1})^{3}u_{h}^{n-1}u_{hxx}^{n}, \,\,\,
\ee_{26}^{n} = (\eta_{h}^{n})^{3}u_{hx}^{n}u_{hx}^{n+1}
- (\eta_{h}^{n-1})^{3}u_{hx}^{n-1}u_{hx}^{n}.
\end{align*}
But
$\ee_{11}^{n} = \eta_{h}^{n}\delta u_{hx}^{n} + u_{hx}^{n}\delta\eta_{h}^{n-1}$, 
$\ee_{12}^{n} = u_{h}^{n}\delta\eta_{hx}^{n} + \eta_{hx}^{n}\delta u_{h}^{n-1}$.
Multiplying the first equation in \eqref{eq447} by $\delta\eta_{h}^{n}$, integrating with respect to
$x$ on $[0,1]$, and noting that 
$(u_{h}^{n}\delta \eta_{hx}^{n},\delta\eta_{h}^{n}) 
= -(u_{hx}^{n}\delta\eta_{h}^{n},\delta\eta_{h}^{n})/2$, we obtain
\[
\tfrac{1}{2}\tfrac{d}{dt}\|\delta\eta_{h}^{n}\|^{2} = -(\eta_{h}^{n}\delta u_{hx}^{n},
\delta\eta_{h}^{n}) +\tfrac{1}{2}(u_{hx}^{n}\delta\eta_{h}^{n},\delta\eta_{h}^{n})
-(u_{hx}^{n}\delta\eta_{h}^{n-1},\delta\eta_{h}^{n})
-(\eta_{hx}^{n}\delta u_{h}^{n-1},\delta\eta_{h}^{n}).
\]
Therefore from the Cauchy-Schwarz inequality and Corollary 4.8, we conclude that there exists a
constant $C=C(M)$ such that for $0\leq t\leq t^{*}$
\begin{equation}
\tfrac{d}{dt}\|\delta\eta_{h}^{n}\|^{2}\leq C(\|\delta\eta_{h}^{n}\|^{2}
+ \|\delta u_{hx}^{n}\|^{2} + \|\delta\eta_{h}^{n-1}\|^{2} + 
\|\delta u_{h}^{n-1}\|^{2}).
\label{eq448}
\end{equation} 
Note now that 
\begin{align*}
\ee_{21}^{n} & =\eta_{h}^{n}\delta u_{ht}^{n} + u_{ht}^{n}\delta\eta_{h}^{n-1},\,\,\,
\ee_{22}^{n} = (\eta_{h}^{n})^{3}\delta u_{htx}^{n} 
+ \omega_{h}^{n}u_{htx}^{n}\delta\eta_{h}^{n-1},\,\,\,
\ee_{23}^{n} = \eta_{h}^{n}\delta\eta_{hx}^{n} + \eta_{hx}^{n}\delta\eta_{h}^{n-1},\\
\ee_{24}^{n} & = \eta_{h}^{n}u_{h}^{n}\delta u_{hx}^{n}
+(\eta_{h}^{n}u_{h}^{n} - \eta_{h}^{n-1}u_{h}^{n-1})u_{hx}^{n}
=\eta_{h}^{n}u_{h}^{n}\delta u_{hx}^{n} + \eta_{h}^{n}u_{hx}^{n}\delta u_{h}^{n-1}
+ u_{h}^{n-1}u_{hx}^{n}\delta\eta_{h}^{n-1},\\
\ee_{25}^{n} & =(\eta_{h}^{n})^{3}u_{h}^{n}\delta u_{hxx}^{n}+
\bigl((\eta_{h}^{n})^{3}u_{h}^{n}-(\eta_{h}^{n-1})^{3}u_{h}^{n-1}\bigr)u_{hxx}^{n}
= (\eta_{h}^{n})^{3}u_{h}^{n}\delta u_{hxx}^{n} 
+ (\eta_{h}^{n})^{3}u_{hxx}^{n}\delta u_{h}^{n-1} \\
& \hspace{255pt} + \omega_{h}^{n}u_{h}^{n-1}u_{hxx}^{n}\delta\eta_{h}^{n-1},\\
\ve_{26}^{n} & = (\eta_{h}^{n})^{3}u_{hx}^{n}\delta u_{hx}^{n}
+ \bigl((\eta_{h}^{n})^{3}u_{hx}^{n} - (\eta_{h}^{n-1})^{3}u_{hx}^{n-1}\bigr)u_{hx}^{n}
= (\eta_{h}^{n})^{3}u_{hx}^{n}\delta u_{hx}^{n}
+ (\eta_{h}^{n})^{3}u_{hx}^{n}\delta u_{hx}^{n-1} \\
& \hspace{255pt} + \omega_{h}^{n}u_{hx}^{n-1}u_{hx}^{n}\delta\eta_{h}^{n-1},
\end{align*}
where we have denoted
$\omega_{h}^{n}:=(\eta_{h}^{n})^{2} + \eta_{h}^{n}\eta_{h}^{n-1}+(\eta_{h}^{n-1})^{2}$. Multiplying now 
the second equation of \eqref{eq447} by $\delta u_{h}^{n}$ and integrating with respect to $x$ on
$[0,1]$, we see that
\begin{align*}
\tfrac{1}{2}\tfrac{d}{dt}\bigl[& (\eta_{h}^{n}\delta u_{h}^{n},\delta u_{h}^{n})
+ \tfrac{1}{3}\bigl((\eta_{h}^{n})^{3}\delta u_{hx}^{n},\delta u_{hx}^{n}\bigr)\bigr]
= \tfrac{1}{2}(\eta_{ht}^{n}\delta u_{h}^{n},\delta u_{h}^{n})+
\tfrac{1}{2}\bigl(\eta_{ht}^{n}(\eta_{h}^{n})^{2}\delta u_{hx}^{n},\delta u_{hx}^{n}\bigr)\\
&-(u_{ht}^{n}\delta\eta_{h}^{n-1},\delta u_{h}^{n})
- \tfrac{1}{3}(\omega_{h}^{n}u_{htx}^{n}\delta\eta_{h}^{n-1},\delta u_{hx}^{n})
-(\eta_{h}^{n}\delta\eta_{hx}^{n}+\eta_{hx}^{n}\delta\eta_{h}^{n-1},\delta u_{h}^{n})\\
& - (\eta_{h}^{n}u_{h}^{n}\delta u_{hx}^{n}+\eta_{h}^{n}u_{hx}^{n}\delta u_{h}^{n-1},
\delta u_{h}^{n})
- (u_{h}^{n-1}u_{hx}^{n}\delta\eta_{h}^{n-1},\delta u_{h}^{n})
- \tfrac{1}{3}\bigl((\eta_{h}^{n})^{3}u_{h}^{n}\delta u_{hxx}^{n},\delta u_{hx}^{n}\bigr)\\
&-\tfrac{1}{3}\bigl((\eta_{h}^{n})^{3}u_{hxx}^{n}\delta u_{h}^{n-1},\delta u_{hx}^{n}\bigr)
-\tfrac{1}{3}(\omega_{h}^{n}u_{h}^{n-1}u_{hxx}^{n}\delta \eta_{h}^{n-1},\delta u_{hx}^{n})\\
& +\tfrac{1}{3}\bigl((\eta_{h}^{n})^{3}u_{hx}^{n}(\delta u_{hx}^{n}+\delta u_{hx}^{n-1}),
\delta u_{hx}^{n}\bigr) 
+ \tfrac{1}{3}(\omega_{h}^{n}u_{hx}^{n-1}u_{hx}^{n}\delta\eta_{h}^{n-1},\delta u_{hx}^{n}).
\end{align*}
Therefore, from the Cauchy-Schwarz inequality, and Corollary 4.8, we conclude that there exists a
constant $C=C(M,M_{1})$ such that for $0\leq t\leq t^{*}$ 
\[
\tfrac{d}{dt}\bigl[(\eta_{h}^{n}\delta u_{h}^{n},\delta u_{h}^{n})
 + \tfrac{1}{3}\bigl((\eta_{h}^{n})^{3}\delta u_{hx}^{n},\delta u_{hx}^{n}\bigr)\bigr]
\leq C(\|\delta \eta_{h}^{n}\|^{2} + \|\delta u_{h}^{n}\|_{1}^{2}
+ \|\delta \eta_{h}^{n-1}\|^{2} + \|\delta u_{h}^{n-1}\|_{1}^{2}).
\]
From this inequality and \eqref{eq448} we get by addition, integration with respect to $t$, and
Lemma 4.1 (noting that $\delta\eta_{h}^{n}(0)=\delta u_{h}^{n}(0)=0$), the inequality
\[
\|\delta\eta_{h}^{n}(t)\|^{2} + \|\delta u_{h}^{n}(t)\|_{1}^{2}
\leq C \int_{0}^{t}(\|\delta\eta_{h}^{n}(\tau)\|^{2} + \|\delta u_{h}^{n}(\tau)\|_{1}^{2}
+ \|\delta\eta_{h}^{n-1}(\tau)\|^{2} + \|\delta u_{h}^{n-1}(\tau)\|_{1}^{2})d\tau,
\]
for $t\leq t^{*}$, and some constant $C=C(M,M_{1},c_{0})$. Consequently 
\begin{equation*}
\begin{aligned}
\|\delta\eta_{h}^{n}(t)\|^{2} + \|\delta u_{h}^{n}(t)\|_{1}^{2}
& \leq Ct\max_{0\leq \tau\leq t}
(\|\delta\eta_{h}^{n}(\tau)\|^{2} + \|\delta u_{h}^{n}(\tau)\|_{1}^{2})\\
&\,\,\,\,\,\, + Ct\max_{0\leq \tau\leq t}(\|\delta\eta_{h}^{n-1}(\tau)\|^{2} 
+ \|\delta u_{h}^{n-1}(\tau)\|_{1}^{2}).
\end{aligned}
\end{equation*}
Hence, choosing $t_{1}^{*}>0$ so that  
$\alpha:=\bigl(Ct_{1}^{*}/(1-Ct_{1}^{*})\bigr)^{1/2} < 1$, (and so that
$t_{1}^{*}\leq t^{*}$), we obtain the conclusion of the Proposition.
\end{proof}
The convergence of $(\eta_{h}^{n},u_{h}^{n})$ to a solution of \eqref{eqsh} for $0\leq t\leq t_{1}^{*}$
follows:
\begin{proposition} Let $\mathcal{U}^{n}_{h}\in \mathcal{S}_{h}\times\mathcal{S}_{h}$ be defined as
$\mathcal{U}^{n}_{h}=(\eta_{h}^{n},u_{h}^{n})$, $n=1,2,\dots$. Then:\\
{\rm{(i)}} The sequence $\{\mathcal{U}^{n}_{h}\}_{n=1}^{\infty}$ is Cauchy with respect to the norm
\[
\vertiii{\mathcal{U}_{h}^{n}}:=\max_{0\leq t\leq t_{1}^{*}}
(\|\eta_{h}^{n}(t)\|^{2} + \|u_{h}^{n}(t)\|^{2}_{1})^{1/2}.
\]
{\rm{(ii)}} The sequence $\{\mathcal{U}^{n}_{h}\}_{n=1}^{\infty}$ converges to a solution 
$\mathcal{U}_{h}:=(\eta_{h},u_{h})$ of the semidiscrete problem \eqref{eqsh}, with respect to the
norm $\vertiii{\cdot}$. This solution satisfies 
$(\eta_{h},u_{h})\in \mathcal{A}_{t_{1}^{*}}$,  $(\eta_{ht},u_{ht})\in \mathcal{B}_{t_{1}^{*}}$, and
$\min_{0\leq x\leq 1}\eta_{h}(x,t)\geq c_{0}/8$, for $0\leq t\leq t_{1}^{*}$
\end{proposition}
\begin{proof} Letting $\delta\mathcal{U}_{h}^{n}:=(\delta\eta_{h}^{n},\delta u_{h}^{n})$, where
$\delta\eta_{h}^{n}$, $\delta u_{h}^{n}$ as in Proposition 4.9, it follows from \eqref{eq446} that
$\sum_{n=1}^{\infty}\vertiii{\delta\mathcal{U}_{h}^{n}}<\infty$, from which (i) follows.\\
(ii) The convergence of $\{\mathcal{U}_{h}^{n}\}_{n=1}^{\infty}$ follows from (i). Define its limit
as $(\eta_{h},u_{h})\in \mathcal{S}_{h}\times \mathcal{S}_{h}$. Since dim$\mathcal{S}_{h}<\infty$
we have that the convergence $\eta_{h}^{n}\to\eta_{h}$, $u_{h}^{n}\to u_{h}$, as $n\to \infty$, is
valid with respect to any norm of $\mathcal{S}_{h}$. Therefore, using the triangle inequality in the
relations $\eta_{h}=\eta_{h}-\eta_{h}^{n} + \eta_{h}^{n}$, $u_{h} = u_{h} - u_{h}^{n} + u_{h}^{n}$,
and taking into account that $(\eta_{h}^{n},u_{h}^{n})\in \mathcal{A}_{t_{1}^{*}}$ for all $n$, 
yields that $(\eta_{h},u_{h})\in \mathcal{A}_{t_{1}^{*}}$. In addition, since 
$\eta_{h}=\eta_{h}-\eta_{h}^{n} + \eta_{h}^{n}$, and since 
$\min_{0\leq x\leq 1}\eta_{h}^{n}(x,t)\geq c_{0}/8$ for $t\leq t_{1}^{*}$ and all $n$, we conclude
that  $\eta_{h}\geq -\|\eta_{h}-\eta_{h}^{n}\|_{\infty} + c_{0}/8$, from which, letting 
$n\to\infty$, it follows that $\min_{0\leq x\leq 1}\eta_{h}(x,t)\geq c_{0}/8$ for all
$t\leq t_{1}^{*}$. Now it holds that $\|\eta_{ht}^{n}-\eta_{ht}\|\to 0$, as $n\to \infty$, uniformly
in $t\in[0,t_{1}^{*}]$. To see this, note that from the first equations of \eqref{eqshn} it follows that $\eta_{ht}^{n}\to -P_{h}(\eta_{h}u_{h})_{x}$, as $n\to \infty$, 
and integrating with respect to $t$ and letting $n\to\infty$, we have that
$\eta_{h}(t)=\eta_{h}(0) - \int_{0}^{t}P_{h}(\eta_{h}u_{h})_{x}(\tau)d\tau$, and therefore
$\eta_{ht}=-P_{h}(\eta_{h}u_{h})_{x}$ for $0\leq t\leq t_{1}^{*}$ (Note in addition that 
$\max_{0\leq \tau\leq t_{1}^{*}}\|\eta_{ht}^{n}(\tau)- \eta_{ht}(\tau)\|\to 0$, $n\to\infty$, and
that $\max_{0\leq \tau\leq t_{1}^{*}}\|\eta_{ht}^{n}(\tau)- \eta_{ht}(\tau)\|_{\infty}\to 0$, since
dim$\mathcal{S}_{h}<\infty$). We conclude that $(\eta_{h},u_{h})$ is a solution of the first
equation of \eqref{eqsh}. In order to show that $u_{ht}^{n}\to u_{ht}$, $n\to \infty$, in $H^{1}$,
uniformly with respect to $t\in [0,t_{1}^{*}]$, define 
\begin{align*}
& b_{h}^{n}:= -P_{h}(\eta_{h}^{n}\eta_{hx}^{n+1}+\eta_{h}^{n}u_{h}^{n}u_{hx}^{n+1})
- F_{h}\bigl((\eta_{h}^{n})^{3}(u_{h}^{n}u_{hxx}^{n}-u_{hx}^{n}u_{hx}^{n+1})\bigr),\\
& b_{h}:=-P_{h}(\eta_{h}\eta_{hx}+\eta_{h}u_{h}u_{hx})
- F_{h}\bigl(\eta_{h}^{3}(u_{h}u_{hxx}-u_{hx}^{2})\bigr),
\end{align*}
and, given $v\in\mathcal{S}_{h}$ such that $\min_{0\leq x\leq 1}v(x)\geq c>0$, define the operator
$A_{h}(v) :\mathcal{S}_{h}\to \mathcal{S}_{h}$ by 
$A_{h}(v)w=P_{h}(vw)+F_{h}(v^{3}w_{x})$, $w\in\mathcal{S}_{h}$. If $f\in\mathcal{S}_{h}$, then by
Lemma 4.1, \eqref{eq42}, $\|A_{h}^{-1}(v)f\|_{1}\leq C\|f\|$, for some constant $C$ depending on $c$
only. Note that the second equation of \eqref{eqshn} is written in the form 
$A_{h}(\eta_{h}^{n})u_{ht}^{n+1} = b_{h}^{n}$. If now $v_{h}^{n}$, $v_{h}$ are the solutions of
the problems $A_{h}(\eta_{h})v_{h}^{n}=b_{h}^{n}$, $A_{h}(\eta_{h})v_{h}=b_{h}$, respectively, then
$\|v_{h}^{n}-v_{h}\|_{1}\leq C_{0}\|b_{h}^{n}-b_{h}\|$, for some constant $C_{0}$ depending on
$c_{0}$. Taking $n\to\infty$ in the last inequality, we see that the sequence 
$\{v_{h}^{n}\}_{n=0}^{\infty}$ converges to $v_{h}$ in $H^{1}$ uniformly in $t\in [0,t_{1}^{*}]$.
Moreover it holds that
\[
A_{h}(\eta_{h})(u_{ht}^{n+1}-v_{h})=[A_{h}(\eta_{h})-A_{h}(\eta_{h}^{n})]u_{ht}^{n+1}
+ b_{h}^{n} - b_{h}.
\]
Hence $u_{ht}^{n+1} - v_{h}=A_{h}^{-1}(\eta_{h})[A_{h}(\eta_{h})-A_{h}(\eta_{h}^{n})]u_{ht}^{n+1}
+ v_{h}^{n} - v_{h}$, and 
\[
\|u_{ht}^{n+1} - v_{h}\|_{1}\leq C_{0}\|[A_{h}(\eta_{h})-A_{h}(\eta_{h}^{n})]u_{ht}^{n+1}\|
+ \|v_{h}^{n} - v_{h}\|_{1}.
\]
From this estimate, since $u_{ht}^{n}\in\mathcal{B}_{t_{1}^{*}}$ for all $n$, we see that
$u_{ht}^{n}\to v_{h}$, as $n\to \infty$, in $H^{1}$, uniformly in $t\in [0,t_{1}^{*}]$. Since now
$u_{h}^{n}-\int_{0}^{t}v_{h}d\tau = \int_{0}^{t}(u_{ht}^{n}-v_{h})d\tau + u_{h}(0)$, taking
$n\to\infty$, we see that $u_{h}=\int_{0}^{t}v_{h}d\tau+u_{h}(0)$, and therefore $u_{ht}=v_{h}$. Thus
$\max_{0\leq t\leq t_{1}^{*}}\|u_{ht}^{n}(\tau)-u_{ht}(\tau)\|_{1}\to 0$, $n\to\infty$. This result
implies that taking $n\to\infty$ in the second equation of \eqref{eqshn}, gives the second equation
of \eqref{eqsh}. Moreover from the convergence of the sequence $\{\eta_{ht}^{n}\}_{n=0}^{\infty}$
that was previously proved, and the fact that dim$\mathcal{S}_{h}<\infty$, we get 
$(\eta_{ht},u_{ht})\in \mathcal{B}_{t_{1}^{*}}$, and the proof of the proposition is now complete.
\end{proof}
We finally establish the uniqueness of the solution $(\eta_{h},u_{h})$ of \eqref{eqsh}, for
$0\leq t\leq t_{1}^{*}$.
\begin{proposition} The semidiscrete problem \eqref{eqsh} has at most one solution 
$(\eta_{h},u_{h})$, for $0\leq t\leq t_{1}^{*}$, such that $(\eta_{h},u_{h})\in\mathcal{A}_{t_{1}^{*}}$, 
$(\eta_{ht},u_{ht})\in\mathcal{B}_{t_{1}^{*}}$,
$\min_{0\leq x\leq 1}\eta_{h}(x,t)\geq c_{0}/8$.
\begin{proof} Let $(\eta_{h}^{a}, u_{h}^{a})$, 
$(\eta_{h}^{b},u_{h}^{b})\in \mathcal{A}_{t_{1}^{*}}$, be two solutions of \eqref{eqsh} such that 
$(\eta_{ht}^{a},u_{ht}^{a})$, $(\eta_{ht}^{b},u_{ht}^{b})\in\mathcal{B}_{t_{1}^{*}}$, with 
$\max_{0\leq x\leq 1}\eta_{h}^{a}(x,t)\geq c_{0}/8$, 
$\max_{0\leq x\leq 1}\eta_{h}^{b}(x,t)\geq c_{0}/8$, for $t\leq t_{1}^{*}$. Putting
\[
\delta\eta_{h}=\eta_{h}^{a}-\eta_{h}^{b}, \quad \delta u_{h}=u_{h}^{a}-u_{h}^{b},
\]
then, from the equations of \eqref{eqsh}, we obtain
\begin{equation}
\begin{aligned}
& \delta\eta_{ht} = -P_{h}(\ve_{1})_{x}, \\
&P_{h}\ve_{2}=-F_{h}(\ve_{3}) -P_{h}(\ve_{4}+\ve_{5}) - F_{h}(\ve_{6}-\ve_{7}),
\end{aligned}
\label{eq449}
\end{equation}
where 
\begin{align*}
\ve_{1} & =\eta_{h}^{a}u_{h}^{a}-\eta_{h}^{b}u_{h}^{b},\hspace{19pt}
\ve_{2} = \eta_{h}^{a}u_{ht}^{a}-\eta_{h}^{b}u_{ht}^{b}, \hspace{36pt}
\ve_{3} = (\eta_{h}^{a})^{3}u_{htx}^{a}-(\eta_{h}^{b})^{3}u_{htx}^{b},\\
\ve_{4} & = \eta_{h}^{a}\eta_{hx}^{a}-\eta_{h}^{b}\eta_{hx}^{b},\,\,\,\,\,\,
\ve_{5} = \eta_{h}^{a}u_{h}^{a}u_{hx}^{a}-\eta_{h}^{b}u_{h}^{b}u_{hx}^{b},\,\,\,\,
\ve_{6} = (\eta_{h}^{a})^{3}u_{h}^{a}u_{hxx}^{a} - (\eta_{h}^{b})^{3}u_{h}^{b}u_{hxx}^{b},\\
\ve_{7} & = (\eta_{h}^{a})^{3}(u_{hx}^{a})^{2} - (\eta_{h}^{b})^{3}(u_{hx}^{b})^{2}.
\end{align*}
Since $\ve_{1} = \eta_{h}^{a}\delta u_{h} + u_{h}^{b}\delta\eta_{h}$, multiplying the first equation
of \eqref{eq449} by $\delta\eta_{h}$, integrating with respect to $x$ on $[0,1]$ and taking into
account that $\bigl((u_{h}^{b}\delta\eta_{h})_{x},\delta\eta_{h}\bigr)
=(u_{hx}^{b}\delta\eta_{h},\delta\eta_{h})/2$, we get 
\[
\tfrac{1}{2}\tfrac{d}{dt}\|\delta\eta_{h}\|^{2}=
-(\eta_{hx}^{a}\delta u_{h},\delta\eta_{hx})-(\eta_{h}^{a}\delta u_{hx},\delta\eta_{h})
-\tfrac{1}{2}(u_{hx}^{b}\delta\eta_{h},\delta\eta_{h}).
\]
Therefore, by the Cauchy-Schwarz inequality and our hypotheses for $\eta_{h}^{a}$, $u_{h}^{b}$, we
conclude that there exists a constant $C(M)$ such that
\begin{equation}
\tfrac{d}{dt}\|\delta\eta_{h}\|^{2}\leq 
C(M)(\|\delta\eta_{h}\|^{2} + \|\delta u_{h}\|_{1}^{2}).
\label{eq450}
\end{equation}
In addition, note that 
\begin{align*}
\ve_{2} & = \eta_{h}^{a}\delta u_{ht} + u_{ht}^{b}\delta\eta_{h},\hspace{94pt}
\ve_{3} = (\eta_{h}^{a})^{3}\delta u_{htx} + \omega_{h}u_{htx}^{b}\delta\eta_{h},\,\,\,
\ve_{4} = \eta_{h}^{a}\delta\eta_{hx} + \eta_{hx}^{b}\delta\eta_{h},\\
\ve_{5} & = \eta_{h}^{a}u_{h}^{a}\delta u_{hx} + \eta_{h}^{a}u_{hx}^{b}\delta u_{h}
+u_{h}^{b}u_{hx}^{b}\delta\eta_{h},\,\,\,
\ve_{6} = (\eta_{h}^{a})^{3}u_{h}^{a}\delta u_{hxx} + 
(\eta_{h}^{b})^{3}u_{hxx}^{b}\delta u_{h} + \omega_{h}u_{h}^{a}u_{hxx}^{b}\delta\eta_{h},\\
\ve_{7} & = (\eta_{h}^{a})^{3}w_{h}\delta u_{hx} + \omega_{h}(u_{hx}^{b})^{2}\delta\eta_{h},
\end{align*}
where $\omega_{h}=(\eta_{h}^{a})^{2} + \eta_{h}^{a}\eta_{h}^{b} + (\eta_{h}^{b})^{2}$, and
$w_{h} = u_{hx}^{a} + u_{hx}^{b}$. Multiplying the second equation of \eqref{eq449} by 
$\delta u_{h}$, integrating with respect to $x$ on $[0,1]$ and using similar estimates on $\ve_{i}$,
$1\leq i\leq 7$ with the ones used e.g. in the course of the proof of Proposition 4.9, we finally
obtain
\begin{equation}
\tfrac{d}{dt}\bigl[(\eta_{h}^{a}\delta u_{h},\delta u_{h}) 
+ \tfrac{1}{3}\bigl((\eta_{h}^{a})^{3}\delta u_{hx},\delta u_{hx}\bigr)\bigr]
\leq C(\|\delta\eta_{h}\|^{2} + \|\delta u_{h}\|_{1}^{2}).
\label{eq451}
\end{equation}
Adding \eqref{eq450} and \eqref{eq451}, integrating with respect to $t$, taking into account the
hypothesis that $\min_{0\leq x\leq 1}\eta_{h}^{a}\geq c_{0}/8$ for $t\leq t_{1}^{*}$, Lemma 4.1,
and that $\delta\eta_{h}(0)=\delta u_{h}(0)=0$, gives
\[
\|\delta\eta_{h}(t)\|^{2} + \|\delta u_{h}(t)\|_{1}^{2} 
\leq C\int_{0}^{t}(\|\delta\eta_{h}(\tau)\|^{2} + \|\delta u_{h}(\tau)\|_{1}^{2})d\tau,
\]
for $t\leq t_{1}^{*}$, and some constant $C=C(M,M_{1},c_{0})$. Therefore, Gronwall's Lemma
implies that $\|\delta\eta_{h}(t)\|=\|\delta u_{h}(t)\|_{1}=0$, giving 
$\eta_{h}^{a}=\eta_{h}^{b}$, $u_{h}^{a}=u_{h}^{b}$ in $[0,1]\times[0,t_{1}^{*}]$.
\end{proof}
\end{proposition}
\bibliographystyle{amsalpha} 

\end{document}